\documentclass[a4paper,12pt]{amsart}
\usepackage{amsmath,
			amssymb,
			amsthm,
			amsfonts,
			mathrsfs, 
			enumerate, 
			mathtools} 
\usepackage[utf8]{inputenc}
\usepackage{graphicx}
\usepackage[all]{xy}
\usepackage[normalem]{ulem}
\usepackage{url} 
\usepackage{color}
\usepackage[top=1.5in,
			bottom=1.5in,
			left=1in,
			right=1in]{geometry}
\usepackage{braket}
\usepackage{geometry}
\usepackage[dvipsnames]{xcolor}

\usepackage{pifont}

\usepackage{mathabx}

\usepackage[all]{xy}

\usepackage{hyperref}
\usepackage{graphicx}

\usepackage{tikz}

\usepackage{txfonts}

\newtheorem{theorem}{Theorem}

\newtheorem{proposition}{Proposition}[section]
\newtheorem{lemma}[proposition]{Lemma}
\newtheorem{corollaire}[proposition]{Corollary}

\theoremstyle{definition}
\newtheorem{definition}[proposition]{Definition}
\newtheorem{example}[proposition]{Example}

\newtheorem{remarque}[proposition]{Remark}
\newtheorem{remarques}[proposition]{Remarks}

\newcommand{\bl}{\begin{lemma}}
\newcommand{\bp}{\begin{proposition}}
\newcommand{\bt}{\begin{theorem}}
\newcommand{\bc}{\begin{corollaire}}
\newcommand{\be}{\begin{equation}}
\newcommand{\bee}{\begin{equation*}}
\newcommand{\bd}{\begin{definition}}
\newcommand{\bdp}{\begin{definitionproposition}}
\newcommand{\bex}{\begin{example}}
\newcommand{\br}{\begin{remarque}}
\newcommand{\bpr}{\begin{proof}}

\newcommand{\el}{\end{lemma}}
\newcommand{\ep}{\end{proposition}}
\newcommand{\et}{\end{theorem}}
\newcommand{\ec}{\end{corollaire}}
\newcommand{\ee}{\end{equation}}
\newcommand{\eee}{\end{equation*}}
\newcommand{\ed}{\end{definition}}
\newcommand{\edp}{\end{definitionproposition}}
\newcommand{\eex}{\end{example}}
\newcommand{\er}{\end{remarque}}
\newcommand{\epr}{\end{proof}}

\newcommand{\secref}[1]{Section~\ref{#1}}
\newcommand{\subsecref}[1]{Subsection~\ref{#1}}
\newcommand{\thmref}[1]{Theorem~\ref{#1}}
\newcommand{\propref}[1]{Proposition~\ref{#1}}
\newcommand{\lemref}[1]{Lemma~\ref{#1}}
\newcommand{\corref}[1]{Corollary~\ref{#1}}

\newcommand{\remref}[1]{Remark~\ref{#1}}

\newcommand{\negro}{\color{black}}

\newcommand\N{\mathbb N}
\newcommand\R{\mathbb R}
\newcommand\C{\mathbb C}
\newcommand\Z{\mathbb Z}

\renewcommand\1{\hbox{\ding{192}}}
\renewcommand\2{\hbox{\ding{193}}}
\newcommand\3{\hbox{\ding{194}}}
\newcommand\4{\hbox{\ding{195}}}

\newcommand{\sbat}{S^1}

\newcommand{\stres}{S^3}
\newcommand{\zdos}{{\Z}_{_2}}
\newcommand{\Msbat}{M^{\sbat}}
\newcommand{\Ssbat}{S^{\sbat}}
\newcommand{\Mstres}{M^{\stres}}

\newcommand{\ib}[2]{{#1}_{_{#2}}}
\newcommand{\Hiru}[3]{{#1}^{^{#2}}_{_{#3}}}
\newcommand{\lau}[4]{{#1}^{^{#2}}_{_{#3}}{\left( #4 \right)}}
\newcommand{\coho}[3]{{#1}^{^{#2}}{\left( #3 \right)}}

\newcommand{\tetra}[4]{{}_{_{#1}} {#2}^{^{#3}}_{_{#4}}       }

\def\per{\overline}
\newcommand{\tq}{\ \big| \ }
\newcommand{\stra}{\mathscr{S}}
\newcommand{\strasing}{\mathscr{S}^{sing}}
\newcommand{\menos}{\smallsetminus}
\newcommand{\rondp}{\raise1pt\hbox{\tiny $\circ$}}
\newcommand\phii{{\raise2pt\hbox{$\varphi$}}}  
  
\newcommand{\wt}{\widetilde} 
\newcommand\BOm{\underline\Omega}
\newcommand\Om{\Omega}  
\newcommand\om{\omega} 
\newcommand{\codim}{\mathop{\rm codim \, } \nolimits} 
\newcommand{\pr}{\mathop{\rm pr \, } \nolimits} 
 
\newcommand{\Ad}{\mathop{\rm Ad \, } \nolimits}
 
\newcommand{\sudos}{\mathfrak{su}(2)}
\newcommand\chii{\raise2pt\hbox{$\chi$}}

\newcommand{\TO}{\longrightarrow}
\newcommand{\Ima}{\mathop{\rm Im \, } \nolimits}
\newcommand\tc{{\mathtt c}}
\newcommand\tv{{\mathtt v}}
\newcommand\tw{{\mathtt w}}
\newcommand\rc{{\mathring{\tc}}}

\usepackage{etoolbox}
\makeatletter
\patchcmd{\@maketitle}
  {\ifx\@empty\@dedicatory}
  {\ifx\@empty\@date \else {\vskip3ex \centering\footnotesize\@date\par\vskip1ex}\fi
   \ifx\@empty\@dedicatory}
  {}{}
\patchcmd{\@adminfootnotes}
  {\ifx\@empty\@date\else \@footnotetext{\@setdate}\fi}
  {}{}{}
\makeatother

\title[Gysin Braid]{The Gysin Braid for 
$\stres$-actions on manifolds}

\date{\today}

\author[J.I.~Royo Prieto]{Jos\'{e} Ignacio Royo Prieto}
\address{Matematika Saila\\ Zientzia eta Teknologia Fakultatea\\ University of the Basque Country UPV/EHU\\ Barrio Sarriena s/n\\ 48940 Leioa\\Spain.}
\email{joseignacio.royo@ehu.eus}
\thanks{Partially supported by Ministerio de Ciencia e Innovaci\'on, Spain, grant PID2022-139631NB-I00. The authors acknowledge that the research cooperation was funded by the program Excellence Initiative Research University at the Jagiellonian University in Krakow within the framework of the research group Reeb-Reinhardt 2022.}

\author{Martintxo Saralegi-Aranguren}
\address{Laboratoire de Math{\'e}matiques de Lens\\  
      EA 2462 \\
      Universit\'e d'Artois\\
         SP18, rue Jean Souvraz\\
          62307 Lens Cedex\\
         France}
\email{martin.saraleguiaranguren@univ-artois.fr}


\begin{document}  

\begin{abstract}
	In a previous work, we constructed a Gysin sequence that relates the cohomology of a manifold $M$ to that of the orbit space $M/\stres$, where the sphere  $\stres$ acts smoothly on $M$. This sequence includes an exotic term that depends on $M^{\sbat}$,
	the subset of points fixed by the action of the subgroup $\sbat$.
	
	The orbit space is a stratified pseudomanifold, which is a type of singular space where intersection cohomology can be applied. When the action is semi-free, the second author has already constructed a Gysin sequence that relates the cohomology of $M$ to the intersection cohomology of $M/\stres$.
	
	However, what happens when the action is not semi-free? This is the main focus of this work. The situation becomes more complex, and we do not find just a Gysin sequence. Instead, we construct a Gysin braid that relates the cohomology of $M$ to the intersection cohomology of $M/\stres$. This braid also contains an exotic term that depends on the intersection cohomology of the fixed point subset $M^{\sbat}$.

\end{abstract}

\maketitle

Given a smooth free action of the sphere $S^3$ on a smooth manifold $M$, we have a sphere bundle and the Gysin sequence
\begin{equation}\label{gys3Clas0}
\xymatrix{
\cdots \ar[r] & \coho H {*-1} M \ar[r] & \coho H {*-4} {M/S^3} \ar[r] & \coho H {*} {M/S^3} \ar[r] & \coho H {*} M \ar[r] &
&\cdots
}
\end{equation}
which relates the cohomologies of the manifold $M$ and that of its orbit space $M/S^3$ (see for example \cite{MR0336651,MR0413144}). In the case of a semi-free action, we do not have a sphere bundle, but we have the Gysin sequence
\begin{equation}\label{gys3Clas}
\xymatrix{
\cdots \ar[r] & \coho H {*-1} M \ar[r] & \coho H {*-4} {M/S^3, M^{\stres}} \ar[r] & \coho H {*} {M/S^3} \ar[r] & \coho H {*} M \ar[r] & \cdots .
&
}
\end{equation}

In the general case, an exotic term appears:
\be\label{NuestraGysin}
\xymatrix@C=7mm{
\cdots \ar[r] & \coho H {*-1} M \ar[r] & \coho H {*-4} {M/S^3, \Sigma/S^3} \oplus \left( \coho H {*-3} {M^{\sbat}}\right)^{-\zdos}\ar[r] & \coho H {} {M/S^3} \ar[r] & \coho H {*} M \ar[r] & \cdots
&
}
\ee
(cf. \cite{MR3119667}). 
In this context, $\Sigma \subset M$ denotes the subset of points in $M$ whose isotropy group is  infinite. The group $\zdos$ acts on $M^{\sbat}$ by $j \in \stres$.

When the action is free, the orbit space is a manifold. In the more general case, it is a pseudomanifold. 
The second author constructed in \cite{MR1184085} the following exact sequence 
\begin{equation}\label{gysfinal}
\xymatrix{
\cdots \ar[r] & \coho H {*-1} M \ar[r]  & \lau H {*-4}{\per p - \per e}  {M/S^3} \ar[r] &  \lau H {*}  {\per p} {M/S^3} \ar[r] & \coho H {*} M \ar[r]  .
&\cdots 
}
\end{equation}
This sequence establishes a connection between the cohomology of the space $M$ and the intersection cohomology of its orbit space $M/\stres$ in the case of a semi-free action. 
This exact sequence holds when  $\per 0 \leq \per p \leq \per t$. The Euler perversity $\per e$ takes the value 4 on singular strata.
 In the special case where $\per p = \per 0$, the exact sequence \eqref{gysfinal} simplifies to \eqref{gys3Clas}.

In this paper, we establish a connection between the cohomology of a manifold $M$ and the intersection cohomology of its orbit space $M/\stres$, for any smooth action $\Phi \colon S^3 \times M \to M$. Of particular interest is the case where the action has three-dimensional orbits, also known as a \emph{mobile action}.

Given a mobile action, we  obtain the following Gysin braid which relates the cohomologies $\coho H*M$ and 
$\lau H {*} {\per p} {M/\stres} $
\bee	
\scalebox{.7}
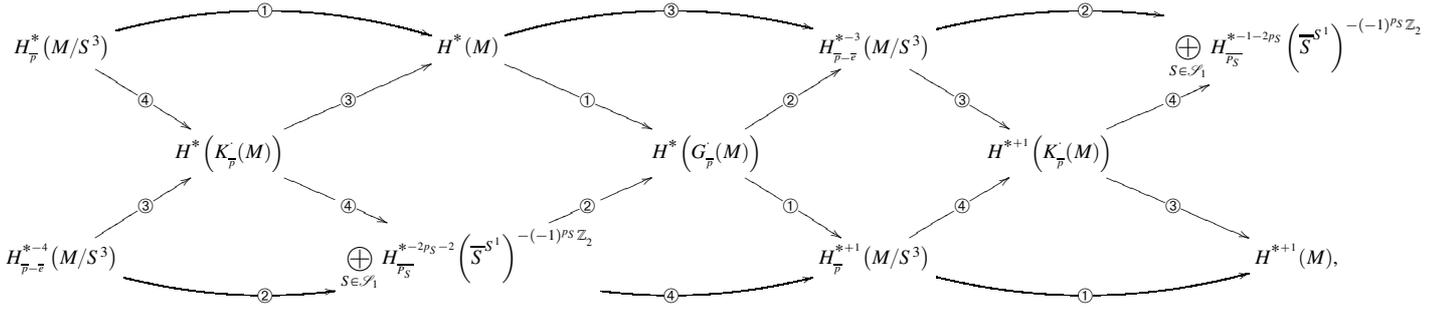
{
 \xymatrix{
\lau H {*} {\per p} {M/\stres} \ar[rd] |\4\ar@/^{7mm}/[rr] |\1 &
&
\coho H*M \ar[rd] |\1\ar@/^{7mm}/[rr] |\3&
&
\lau H {*-3} {\per p - \per e} {M/\stres} \ar[rd] |\3 \ar@/^{7mm}/[rr] |\2&
&
\displaystyle \bigoplus_{S \in \stra_1} \lau H {*-1-2p_S} {\per{P_S}} {\per S^{\sbat}}^{-(-1)^{p_S}\zdos}  \\
&
\coho H * { \lau K \cdot {\per p} M}  \ar[ur]  |\3\ar[dr] |\4 &
&
\coho H {*} { \lau {G} \cdot {\per p } M } \ar[ur]  |\2 \ar[dr] |\1&
&
\coho H {*+1} { \lau K \cdot {\per p} M}   \ar[ur] |\4\ar[dr] |\3 &
\\
\lau H {*-4} {\per p - \per e} {M/\stres}  \ar[ur] |\3\ar@/_{7mm}/[rr] |\2&
&
\displaystyle \bigoplus_{S \in \stra_1} \lau H {*-2p_S-2} {\per{P_S}} {\per S^{\sbat}}^{-(-1)^{p_S}\zdos} \ar[ur] |\2\ar@/_{7mm}/[rr]  |\4&
&
\lau H {*+1}  {\per p} {M/\stres}\ar[ur] |\4\ar@/_{7mm}/[rr] |\1&
&
\coho H{*+1}M,
}}
\eee
where
 \begin{itemize}
 
 \item $\stra_1$ (resp. $\stra_3$) is the family of singular strata $S$ of $M$ with $\dim \stres_x =1$ (resp. 3) for any $x \in S$.
 Here $\stres_x$ denotes the isotropy subgroup.
 
 \item the Euler perversity $\per e$ takes the value 4 (resp. 2) on the strata of $\stra_1$ (resp $\stra_3$),
 
 \item the perversity $\per p$ lies between $\per 0$ and the top perversity $\per t$ on
 $M$,
 
 \item the number $p_S$ denotes the integer part  of $\per{p}(S)/2 $,
 
 \item the perversity $\per{P_S}$ on the filtered space $\per S^{\sbat}$ is defined by $\per{P_S} (Q) = \per p(Q) - 2 p_S - 2$ for any $Q \in \stra_3$ with $Q \subset \per S$,
 
 \item the  Gysin term  $\lau G * {\per p} M$ represents the cokernel of the map induced by the natural 	 projection $\pi \colon M \to M/S^3$ and
 
 \item   the co-Gysin term $\lau K * {\per p} M$
denotes the kernel of the map induced by integrating along the fibers of $\pi$
 
 \end{itemize}
 (see \thmref{C}).  The braid consists of four long exact sequences, denoted by $\1$, $\2$, $\3$, and $\4$. All the triangles and diamonds in the braid are commutative.
The top and bottom sequences in the braid are semi-exact  and both have the same exactness defaults (cf. \remref{Rem} (c)).
The cohomologies of the Gysin and co-Gysin terms are interconnected through  the long exact sequences of \remref{Rem} (b).

 \bigskip
 
Let's analyze the four Gysin sequences that make up the  Gysin braid.
 
 \smallskip
 
 In the classical framework of a free action, there exist two methods for constructing the Gysin sequence \eqref{gys3Clas0}. One approach involves employing the pullback induced by the natural projection $\pi \colon M \to M/S^3$, while the other entails integrating along the fibers of $\pi$. Remarkably, both methodologies yield identical outcomes, resulting in the Gysin sequence.
Meanwhile, the Gysin sequence \eqref{NuestraGysin} is derived via the former method.

 \medskip

In the broader context of this paper's discussion on mobile actions, we utilize both techniques, and they yield distinct results. This fundamental difference is the primary reason why the Gysin braid appears instead of a Gysin sequence.

\smallskip

\begin{itemize}
 \item[$\bullet$] The pullback   associated to the projection $\pi$  induces the long exact sequence $\1$:
 $$
\xymatrix{
\cdots \ar[r] & \coho H {*-1} M \ar[r] & \coho H {*-1} {\lau G \cdot {\per p} M}\ar[r] & \lau H {*} {\per p} {M/S^3} \ar[r] & \coho H {*} M \ar[r]
&\cdots
}
$$
where the \emph{Gysin term} $\lau {G} * {\per p} M$ is the cokernel of $\pi^*$.
This is the first Gysin sequence associated to the action $\Phi$.

We can determine the cohomology of the Gysin term through the sequence $\2$, which employs integration along the fibers of $\pi$. 
This method is used in \cite{MR3119667}, where we implicitly work with the perversity $\per p = \per 0$. In this context, we obtain the Gysin sequence \eqref{NuestraGysin} since the sequence $\2$ splits at the position of the connecting map. However, it's important to note that for other perversities, the sequence $\2$ may not necessarily split as demonstrated by the example in \remref{Rem} (d).

\item[]

 \item[$\bullet$] Employing the integration along the fibers of $\pi$, we obtain the long exact sequence $\3$:
$$
\xymatrix{
\cdots \ar[r] & \coho H {*-1} M \ar[r] & \lau H {*-4} {\per p - \per e} {M/S^3} \ar[r] & \coho H {*} {\lau K \cdot {\per p} M}\ \ar[r] &
\coho H {*} M \ar[r]
&\cdots,
}
$$
where the \emph{co-Gysin term} $\lau {K} * {\per p} M$ is the kernel of the integration operator along the fibers of $\pi$.
This is the second Gysin sequence associated to the action $\Phi$.

 We compute the cohomology of the co-Gysin term using the sequence $\4$, which employs the pullback operator induced by $\pi$. It's important to note that unlike the previous case, sequence $\3$ does not necessarily split, even for the perversity $\per p=\per 0$, as demonstrated in the example in \remref{Rem} (d).

\item[]

 \item[$\bullet$] 
If the exotic term $\displaystyle \bigoplus_{S \in \stra_1} \lau H {*} {\per{P_S}} {\per S^{\sbat}}^{-(-1)^{p_S}\zdos}$ vanishes, then the Gysin braid simplifies to the sequence \eqref{gysfinal}. In particular, this happens  when the action $\Phi$ is semi-free, as noted in \remref{Rem} (a).

\item[]

 \item[$\bullet$] 
  Another approach to constructing a Gysin sequence involves utilizing the Leray-de Rham spectral sequence. Recall that in the case of a differentiable action $\Phi \colon G \times M \to M$ of a connected compact Lie group $G$ on a manifold $M$, there exists a spectral sequence $\Hiru E {i,j} r $ converging to $\coho H {i+j} M$, where $\Hiru E {i,j} 2 = \coho H i{M/G} \otimes \coho HjG$.
When $G =S^3$ and the action is free, almost-free, or semifree, this spectral sequence degenerates into a Gysin sequence (cases \eqref{gys3Clas0}, \eqref{gys3Clas}, and \eqref{gysfinal}). This is because the second term of the spectral sequence contains only two levels: $j=0,3$. However, the situation becomes significantly more complex when the action is not free, and computing the second term of this spectral sequence becomes challenging (see \cite{MR1621397}).

In the case of a mobile action of $S^3$, and using singular cohomology, we have shown in \cite{MR3119667} that the second term of the spectral sequence possesses  three levels ($j=0,2,3$) and that the spectral sequence degenerates into the Gysin sequence \eqref{NuestraGysin}.
 In this paper, we prove that this phenomenon persists within the framework of intersection cohomology. The Leray-de Rham spectral sequence 
 $\tetra {\per p}E {i,j} r $ depends on a perversity $\per 0 \leq \per p \leq \per t$ on $M$ and satisfies the following properties.

 \begin{itemize}
 \item[$\star$] It converges: $\tetra {\per p}E {i,j} r\Rightarrow H^{i+j}(M)$.
 \item[$\star$] The second page is given by
 $$
 \tetra {\per p}E {i,j} 2 =
 \left\{
 \begin{array}{ll}
 \lau H i {\per p} {M/\stres} & \hbox{if } j=0\\[,2cm]
  \displaystyle \bigoplus_{S \in \stra_1} \lau H {i-2p_S} {\per{P_S}} {\per S^{\sbat}}^{-(-1)^{p_S}\zdos} & \hbox{if } j=2\\[,6cm]
  \lau H {i} {\per p - \per e} {M/\stres}& \hbox{if } j=3.\\
 \end{array}
 \right.
 $$
 It is $0$ otherwise.

 \item[$\star$] The Gysin term appears in this spectral sequence through  the long exact sequence $\2$ 
 $$
\xymatrix@C=7mm{
\cdots \ar[r] & \coho H {i-1} { \lau {G} \cdot {\per p } M } \ar[r]  &
 \tetra {\per p} E {i-4,3} 2
   \ar[r]^{d_2} &  
 \tetra {\per p} E {i-2,2} 2
\ar[r] &  \coho H {i} { \lau {G} \cdot {\per p } M } \ar[r]  & \cdots .
&
}
$$

\item[$\star$] This spectral sequence degenerates at the third page and produces the long exact sequence $\1$.

 \end{itemize}
 
 \bigskip

\noindent In other words, the information in the Leray-de Rham spectral sequence beyond page $ \tetra {\per p} E {} 2$ is contained within the Gysin braid.

\item[]

 \item[$\bullet$] 
 
The non-mobile actions are simpler, and we obtain
$
\lau H * {} M = \lau H * {\per p}  {M/\stres}
\oplus
 \lau H {*-2} {\per p - \per 2}  {M^{\sbat} }^{-\zdos}$
 (see \secref{non-mobile}).
 
\end{itemize}

\bigskip

We denote by $S^3$ the unit quaternionic group, i.e., $S^3 = \{x_0+x_1i+x_2j + x_3 k \tq x_0^2+x_1^2+x_2^2+x_3^2=1\}$. The distinguished copy of the unit complex numbers is $S^1 =  \{x_0+x_1i \tq x_0^2+x_1^2=1\}$. Notice that the normalizer $N(S^1)$ of $S^1$ in $S^3$ is $\sbat \sqcup j\ \sbat \cong  O(2)$.

 In the following, we consider a smooth action $\Phi \colon S^3 \times M \to M$, where $M$ is a second countable, Hausdorff, smooth manifold of dimension $m$ without boundary.  For the definitions and properties related to compact Lie group actions, we refer the reader to \cite{MR0413144}.

The first section of this work is devoted to studying the intersection cohomology of the orbit space $M/S^3$. We demonstrate how to compute this cohomology using differential forms defined on an open subset of $M$. The complex of invariant intersection forms of $M$ is a key tool for constructing the Gysin braid, which is discussed in the second section. The final two sections of this work focus on constructing and analyzing the Gysin braid associated to the action. This braid arises from integrating along the orbits of the action discussed in section three.

The authors want to express their deep gratitude to the referee for so many remarks, corrections and suggestions, which helped to improve the quality of the paper.

\tableofcontents

\section{Intersection cohomology}

The intersection cohomology of the orbit space $M/S^3$ is originally defined using singular simplices. In this section, we demonstrate an alternative method for computing this cohomology by utilizing differential forms defined on the regular part of $M$.

\subsection{Filtered spaces \cite{CST1}.}\label{FFSS}
Let us consider the 
 partition of $M$ obtained by the following  equivalence relation:
 $$x \ \sim \ y \Leftrightarrow \dim S^3_x =\dim S^3_y.$$
This condition is equivalent to  $(S^3_x)_0$ and $(S^3_y)_0$ being conjugated, where $(-)_0$ denotes the connected component containing the unity. 

 A {\em stratum} is a connected component of any element of this partition. The {\em orbit type stratification} ${\mathscr S}$ of $M$  is the family of these strata.

There are four possible isotropy subgroups of a point in $M$, up to conjugacy: a finite subgroup of $\stres$, $S^1$, the normalizer $N=N(S^1)$ of $S^1$ in $S^3$, and $S^3$ itself  (cf. \cite[Th.~ 9.5, p~.153]{MR0413144}). 

We define ${\mathscr S}={\mathscr S}_0\sqcup {\mathscr S}_1\sqcup {\mathscr S}_3$ as follows:
\bee
\begin{array}{lll}
{\mathscr S}_0= \{ S \in  {\mathscr S}\tq  \dim \stres_x  =0, x \in S \}  & {\mathscr S}_1= \{ S \in  {\mathscr S}\tq \dim \stres_x =1, \  x \in S\} &   {\mathscr S}_3= \{ S \in  {\mathscr S}\tq  \dim \stres_x = 3, \  x \in S\} \\
 \hfill  \hbox{\small \em Mobile strata} \hfill & \hfill  \hbox{\small \em Semi-mobile strata} \hfill &   \hfill  \hbox{\small \em Fixed strata}\hfill 
\end{array}
\eee
The action is considered {\em mobile} if ${\mathscr S}_0\neq \emptyset$. If ${\mathscr S}_0=\emptyset$ and ${\mathscr S}_1\neq\emptyset$, we say that the action is {\em semi-mobile}. The remaining case is the {\em trivial} action.
The set of {\em singular strata} is denoted by $\strasing$. If the action is mobile, then $\strasing=\stra_1\sqcup \stra_3$. If the action is semi-mobile, then $\strasing=\stra_3$.
The union of the singular strata of $\stra$ is  the {\em singular part}  $\Sigma$ of $M$.
Its complementary $M\menos \Sigma$ is the {\em regular part }of $M$.

We define
$
 F_3 = \sqcup_{S \in \stra_3} S = \Mstres $  and $F_1 = \sqcup_{S \in \stra_1} S,
 $
 which are $\stres$-invariant proper submanifolds\footnote{The next Proposition shows that they are manifolds. In fact, these manifolds may have  connected
components with different dimensions.} of $M$. Note that $F_1$ is actually the twisted product $S^3 \times_N (\Msbat\smallsetminus \Mstres)$, where $S^3$ acts on the left of the left factor (see \cite[Th.~5.9]{MR0413144}). Furthermore, if $S \in \stra_1$ we have
\bee\label{eq:eseuno}
S = S^3 \times_N \Ssbat.
\eee
  The union of singular strata is $\Sigma = F_1 \sqcup F_3$ (resp. $F_3$) when the action is mobile (resp. semi-mobile). 
  
 Filtered spaces provide the essential framework for defining singular intersection cohomology, which is the dual of the intersection homology introduced in \cite{MR572580} (see, for example, \cite{CST5}).

\begin{proposition}\label{pro:straM}
The strata of 	 $\stra$ are invariant proper submanifolds.
For each integer $i$ we define \\ $M_i = \sqcup \left\{ S \in \stra \tq \dim S \leq i \right\}$. The filtration
$$\emptyset = M_{-1} \subseteq  M_0 \subseteq \cdots \subseteq  M_i \subseteq \cdots \subseteq M_{m} = M$$
defines a filtered space.

For each integer $i$ we define 
$\left(M/S^3\right)_i = \sqcup \left\{ \pi(S) \tq S \in \stra \hbox{ and } \dim \pi(S) \leq i\right\}$, 
where $\pi \colon M\to M/S^3$ denotes the canonical projection. The filtration
$$\emptyset = \left(M/S^3\right)_{-1} \subseteq  \left(M/S^3\right)_0 \subseteq  \cdots \subseteq \left(M/S^3\right)_i \subseteq \cdots \subseteq \left(M/S^3\right)_{n} = M/S^3$$
defines a filtered structure in $M/S^3$. \end{proposition}
\begin{proof}
Let $S$ be a stratum of $\stra$. Each point $x \in S$ possesses an open neighborhood $S^3$-equivariantely diffeomorphic to the twisted product $S^3 \times_H \R^a$ where the isotropy subgroup 
$H=S^3_x$ acts orthogonally on $\R^a$. The point $x$ becomes the class $\langle1,0\rangle$. Recall that the isotropy subgroup of a point $\langle g,u\rangle \in S^3 \times_H \R^a$ is $gH_u g^{-1}$. So, the trace of $S$ in this neighborhood is $S^3\times_H \R^b$, where $\R^b = \{ u \in \R^a \tq \dim H_u = \dim H\} = \{ u \in \R^a \tq H_0 \cdot u = u\}$. The stratum $S$ is an invariant proper submanifold with $\dim S = 3+b-\dim H$.

 It remains to prove that each $M_i$ and $\left(M/S^3\right)_i$ are closed subsets.  It suffices to verify that  the maps $\dim \colon M \to \Z$ and  $i \colon M \to \Z$, defined by $\dim (x) = \dim S$ and  $i(x) = \dim \pi(S) = \dim S/S^3 =\dim S + \dim S^3_x - 3$, with $S\in \stra$ and $x \in S$, are lower semi-continuous. Since the problem is a local question then we can suppose that $M$ is $S^3\times_H \R^a$. We prove that the functions $\dim$ and $i$ are bigger than $\dim(x)$ and $i(x)$ respectively.
Notice that the map $\langle g,u \rangle \mapsto - \dim S^3_{\langle g,u \rangle}$ is a lower semi-continuous  map since $-\dim S^3_{\langle g,u \rangle} = -\dim H_u \geq -\dim H = -\dim S^3_x$. So, it remains to   study the function $\dim$.
 
  Considering the $G$-equivariant covering $S^3\times_{H_0} \R^a \to S^3\times_H \R^a$, we can suppose that $H$ is connected. Let $\R^a = \R^b \times \R^c$ be the $H$-equivariant orthogonal decomposition of $\R^a$.
 This gives $S^3\times_H \R^a = \R^b \times \left(S^3\times_H \R^c\right)$.
   Given a point $y=\langle g,u \rangle \in S^3 \times_H \R^a$  we consider $Q \in \stra$ the stratum containing this point. In fact, $Q =\R^b \times \left(S^3 \times_H \R^d\right)$ where
$\R^d = \{ v \in \R^a \tq \dim H_v = \dim H_u \} $. We have finished since $\dim (y) =\dim Q =  b+3+d-\dim H \geq b+3-\dim H = \dim S = \dim(x)$.
\end{proof}

 The dimension $m$ of the filtered space $M$ is $\dim M$. The dimension $n$ of the filtered space $M/\stres$ is $m-3$ (resp. $m-1$) when the action is mobile (resp. semi-mobile).

Brylinski-Goresky-MacPherson  showed how to compute the intersection cohomology with differential forms (cf. \cite{MR1197827}). To this effect, they use the Thom-Mather systems.

\subsection{Thom-Mather systems} \label{TMs}
Since $F_1$ and $F_3$  are $\stres$-invariant proper submanifolds of $M$,
we can consider
$\tau_k\colon T_k\to F_k$ 
two $\stres$-invariant  tubular neighborhoods of $F_k$ in $M$, $k=1,3$. 
  Associated
to these tubular neighborhoods we have the following maps:

\medskip

$\bullet$ The {\em radius map} $\rho_k \colon T_k \to [0,\infty)$ is defined fiberwise by $u \mapsto \|u\|$. This map is invariant and smooth.

\medskip

$\bullet$
The {\em dilatation map} $\partial_k \colon [0,\infty) \times T_k \to T_k$,
defined fiberwise by $(t,u) \mapsto t \cdot u$.  It is a smooth  equivariant map.

Given $S \in \stra$ contained in $F_k$ for $k=1,3$, we can define $T_S = \tau_k^{-1}(S)$ and $\tau_S \colon T_S \to S$ as the restriction of $\tau_k$. We can define the maps $\rho_S$ and $\partial_S$ analogously. The {\em soul} of $T_S$ is defined as the open subset $D_S = \rho_S^{-1}([0,2))$.

The family of tubular neighborhoods $\mathfrak{T}_M = \set{T_1,T_3}$ is called a {\em Thom-Mather system of $M$} when

\begin{equation}\label{TMe}
  \left\{
   \begin{array}{c}    \tau_{3} = \tau_{3} \rondp \tau_{1}\\   \rho_{3} = \rho_{3} \rondp \tau_{1}  \end{array} \right\} 
   \hbox{ on }
   T_{1} \cap T_{3} = \tau_{1}^{-1}(T_{3} \cap
   F_1).
   \end{equation}
We have proved in \cite{MR3119667} that there exists an $S^3$-invariant Thom-Mather system of $M$.

Consider  the induced maps $\wt {\tau_k} \colon T_k/\stres \to F_k/\stres$, $\wt{\rho_k} \colon \wt{T_k} \to [0,\infty)$
and 
$\wt {\partial}_k \colon [0,\infty) \times T_k/\stres \to T_k/\stres$.
The family of tubular neighborhoods $\mathfrak{T} _{M/\stres}= \{T_1/\stres,T_3/\stres \}$ is 
a {\em Thom-Mather system of $M/\stres$}.

We need a more precise description of the atlas of the bundle $\tau_3$. The open tubular neighborhood $T_3$ can be chosen as a disjoint union $T_3 = \sqcup \{T_S \tq S \in \mathscr S_3\}$ with $T_S \cap T_{S'} =\emptyset$ if $S \ne S'$. 
There exists an $\stres$-equivariant atlas $\mathcal A = \{ \phii \colon \tau_S^{-1} (U)\to U \times \R^{b+1}\}$ relatively to an orthogonal action $\Phi_S \colon \stres \times \R^{b+1} \to \R^{b+1}$, having the origin as the only fixed point.

\subsection{Intersection differential forms on $M$} \label{13}
The {\em perverse degree} of a differential form $\om \in\coho{\Om}{*}{M\menos\Sigma}$, defined on the regular part $M\menos\Sigma$ of $M$,  relatively to a singular  stratum $S \in \mathscr{S}$ (see \subsecref{FFSS}),
is the number
\bee
||\om||_{S} = \min \left\{ \ell \in \N \tq \om (v_0, \ldots , v_\ell,-) = 0  \ \ 
\hbox{where $v_0, \ldots, v_\ell$ are vectors tangent to the fibers of $ \tau_S \colon D_S \to S$} 
\right\}
\eee
if $\om\ne 0$ on $D_S$. If $\omega=0$ on $D_S$ we define $||\omega||_S = -\infty$. 
The condition $||\omega||_S =||d \omega||_S=0$ is equivalent to stating that the restriction of $\omega$ to $D_S$ is a $\tau_S$-basic form.

We shall need the following properties of the perverse degree:
\be\label{PropGradPer}
\begin{array}{lll}
||\om||_S \leq |\om|, \hbox{degree of $\omega$.}  &
||d\om||_S \leq ||\om||_S +1. &||g^*\omega||_S = ||\omega||_S \ \ \hbox{ for each } g \in \stres. \\[,2cm]
||\om + \eta||_S \leq \max (||\om||_S , ||\eta||_S).&
||\om \wedge \eta||_S \leq  ||\om||_S + ||\eta||_S. & ||i_X \omega||_S \leq ||\omega||_S \hbox{ for each vector field } X.
\\[,2cm]
\end{array}
\ee
These properties come directly from the definition of perverse degree. We have also used   the fact that the Thom-Mather system is $\stres$-invariant.

A {\em perversity} is a map $\per{p} \colon \strasing \to \per \Z =  \Z \sqcup \{-\infty,\infty\}$. The {\em constant perversity} is   $\per{\ell}(S)=\ell$, with $\ell \in \per  \Z $, for any singular stratum $S$.
The {\em top perversity} is defined by $\per{t}(S)=\codim  S -2$ for any singular stratum $S$.

The complex of {\em$\per p$- intersection differential forms} of $M$, relatively to the perversity $\per{p}$
  is defined by
\begin{equation*}
\lau{\Om}{*}{\per{p}}{M}=\left\{ \om \in \coho{\Om}{*}{M\menos\Sigma } \tq \max (||\om||_S,||d\om||_S) \leq \per{p}(S) \ \forall S \in {\mathscr S}
\right\}.
\end{equation*}

The complex $\lau{\Om}{*}{\per{p}}{M}$ computes the cohomology $\lau{H}{*}{}{M}$ for Goresky-MacPherson perversities  \cite{MR1197827}  or for perversities verifying  $\per 0 \leq \per p\leq \per t$ \cite{MR1143404,MR2210257}. The cohomology of this complex is $\coho H * {M, \Sigma}$ (resp. $\coho H * {M
\menos\Sigma}$) when $\per p < \per 0$ (resp. $\per p > \per t$).

We can compute the intersection cohomology of the orbit space $M/S^3$ by using differential forms defined on $M\smallsetminus\Sigma$. To see how, first notice that the natural projection $\pi$ establishes a bijection between the strata of $M$ and those of $M/S^3$. Therefore, a perversity $\per p$ on $M/S^3$ determines a perversity $\per p$ on $M$ as well, which we will still denote by $\per p$, following the formula $\per p (S) = \per p (S/S^3)$. The reverse is also true. Throughout this work, we will consider all perversities on $M$. Let us set such a perversity as $\per p$. 

A {\em basic differential form} is a 
 differential form $\omega \in  \coho \Om*{M\menos \Sigma}$ verifying the following condition: 
 $$
 \omega(v,-)= d\omega(v,-)=0 
 \hbox{ for each vector  $v$ tangent to the fibers of $\pi \colon M\menos \Sigma \to (M\menos \Sigma)/S^3$}.
 $$
The {\em complex of basic differential forms}, denoted by $\lau \Om * {\rm bas} {M\menos \Sigma}$, plays  the role of the complex of differential forms on the orbit space $(M\menos \Sigma)/\stres$, but its forms are defined on the manifold $M\menos \Sigma$. The cohomology of the complex of {\em basic ${\per p}$-intersection forms} $\lau \Om * {\per p} {M/S^3} = \lau \Om * {\rm bas} {M\menos \Sigma} \cap \lau \Om * {\per p}  M$ is denoted by $\lau H * {\per p} {M/S^3}$.
It has been proved in \cite{SW-BICConFib} that this cohomology is isomorphic to  the  $D\per p$-intersection cohomology of the pseudomanifold $M/\stres$.
Notice that this cohomology does not depend on the Thom-Mather system we have chosen.

 Given two perversities $\per q \leq \per p $ on $M$ the {\em step complex} $\lau \Om * {\per p /  \per q} {M/\stres}$ is the quotient $\lau\Om * {\per p } {M/\stres}/ \lau\Om * {\per q} {M/\stres}$ and its cohomology is denoted by 
  $\lau H * {\per p /  \per q} {M/\stres}$ (cf. \cite{MR642001}). 
  
\subsection{Closure of a stratum}\label{ss1}
The exotic term of the Gysin braid we construct in this work uses a particular filtered space we describe now.

Consider a non-closed stratum $S \in \stra_1$. The closure of $S$  is the union of $S$ itself and a family of fixed strata  $\{Q_\ell \in \stra_3\tq \ell \in J\}$, i.e., 
$\per S = S \sqcup_{\ell \in J} Q_\ell$. It is a filtered space whose regular stratum is $S$. Any perversity $\per p$ on $M$ induces a perversity on $\per S$, still denoted by $\per p$, which is defined by the numbers $\per p(Q_\ell)$, $\ell\in J$.
The Thom-Mather system $\mathfrak T_M$ of $M$ induces the Thom-Mather system $\mathfrak T_{\per S} =
\{ T_{Q_\ell} \cap \per S\tq \ell \in J\}$.

In fact, we need to go a step further and consider the space $\per S^{\sbat}$ which is the union $S^{\sbat} \sqcup_{\ell \in J} Q_\ell$.  It is a filtered space whose regular part is $S^{\sbat}$. Any perversity $\per p$ on $M$ induces a perversity on $\per S^{\sbat}$, still denoted by $\per p$, which is defined by the numbers $\per p(Q_\ell)$, $\ell\in J$.
The Thom-Mather system $\mathfrak T_M$ of $M$ induces the Thom-Mather system $\mathfrak T_{\per S^{\sbat}} =
\left\{ T_{Q_\ell} \cap \per S^{\sbat}\tq \ell \in J \right\}$.

The complex of {\em intersection differential forms}   of $\per S^{\sbat}$, relatively to the perversity $\per{p}$, can be defined as in the previous section.
It computes the intersection cohomology $\lau H * {\per p} {\per S^{\sbat}}$. Since $j^2=-1\in \sbat$, the group  $\zdos$ acts on $\per S^{\sbat}$ by $g\cdot x =j(x)$, where $g$ denotes the generator of $\zdos$. Then the group $\zdos$ also acts on
 this cohomology. We shall use the notation
$$\lau H*{\per p} {\per S^{\sbat}}^{-\zdos} = 
\left\{ \omega \in \lau H*{\per p} {\per S^{\sbat}} \tq g \cdot \omega =-\omega \right\}.
$$

In fact, we are going to use a particular perversity on this space.
Associated to any perversity $ \per p $ on $M$ we have the perversity $\per{P_S}$ on 
$\per S^{\sbat}$ defined by
$$
\per{P_S} (Q) = 
\left\{ 
\begin{array}{cl}
\per p(Q) - 2 p_S - 2& \hbox{if $Q = Q_\ell$ for some $\ell \in J$}.\\
0 & \hbox{otherwise}
\end{array}
\right.
$$
where $p_S$ is  the integer part of $\per p(S)/2$.

\subsection{$\sbat$-actions}\label{circulo}
A similar study can be done for a smooth action $\Psi \colon \sbat \times M \to M$. In this case the orbit type stratification is  ${\mathscr S} = {\mathscr S}_0 \sqcup {\mathscr S}_1  $ where
\bee
\begin{array}{lll}
{\mathscr S}_0= \{ S \in  {\mathscr S}\tq  \dim S^1_x  =0, x \in S \}  &\text{and} & {\mathscr S}_1= \{ S \in  {\mathscr S}\tq \dim \sbat_x =1, \  x \in S\}.  \\
 \hfill  \hbox{\small \em Mobile strata (regular strata)} \hfill & &   \hfill  \hbox{\small \em Fixed strata  (singular strata)}\hfill 
\end{array}
\eee
We suppose that the action is mobile, meaning $\stra_0\ne\emptyset$ or non trivial.

The family of singular strata is denoted by $\strasing$ and its union is $\Sigma$.
The manifold $M$ and the orbit space $M/\sbat$ are filtered spaces. Thom-Mather systems also exist. In this case $\mathfrak T_M = \{ T_S \tq S \in \stra_1\}$ where each $\tau_S \colon T_S \to S$ is an $\sbat$-fiber bundle.
The elements of $\mathfrak T_M $ can be chosen to be disjoint.

We need a more precise description of the atlas of the bundle $\tau_S$. There exists an $\sbat$-equivariant atlas $\mathcal A = \{ \phii \colon \tau_S^{-1} (U)\to U \times \R^{2b+2}\}$ relatively to an orthogonal action $\Phi_S \colon \sbat \times \R^{2b+2} \to \R^{2b}+2$, $b\geq 0$, having the origin as the only fixed point.

Proceeding as in \subsecref{13}, we define the complex $\lau \Om*{\per p} M$, which computes the cohomology of $M$. We also define the complex $\lau \Om*{\per p} {M/\sbat}$, whose cohomology $\lau H*{\per p} {M/\sbat}$ is isomorphic to the  $D\per p$-intersection cohomology  of the pseudomanifold $ M/\sbat$.
Notice that this cohomology does not depend on the chosen Thom-Mather system.

 Given two perversities $\per q \leq \per p $ on $M$ the {\em step complex} $\lau \Om * {\per p /  \per q} {M/\sbat}$ is the quotient $\lau\Om * {\per p } {M/\sbat}/ \lau\Om * {\per q} {M/\sbat}$ and its cohomology is denoted by 
  $\lau H * {\per p /  \per q} {M/\sbat}$ (cf. \cite{MR642001}). This cohomology fits into the long exact sequence
  \begin{equation}\label{pasito}
  \xymatrix{
  \cdots \ar[r] & \lau H \ell {\per q}{M/\sbat} \ar[r] & \lau H \ell {\per p}{M/\sbat} \ar[r] & 
   \lau H \ell {\per p/\per q}{M/\sbat}  \ar[r] & \lau H {\ell+1} {\per q}{M/\sbat} \ar[r]& \cdots .
   }
  \end{equation}

\subsection{$N$-actions}\label{circuloN} A  smooth action $\Theta \colon N \times M \to M$ 
induces a circle action $\Psi \colon \sbat \times M \to M$. Since the stratification $\stra_1$ of this last action is $N$-invariant then we can choose an $N$-invariant  Thom-Mather system $\mathfrak T_M = \{ T_{S} \tq S \in \stra_1 \}$, that is, the map $g \colon T_{S} \to T_{g(S)}$ is an $\sbat$-morphism bundle for each $g \in N$.

Since $j^2 = -1\in \sbat$,  the singular part of the $\sbat$-action $\Sigma = \sqcup \{ S \in \stra_1\}$ is $\zdos$-invariant relatively to the action $g\cdot x =j(x)$. Also,  the union $S \cup j(S)$ is $\zdos$-invariant for any $S \in \stra_1$. Notice that we have two possibilities: either $j(S)=S$ or $j(S)\cap S=\emptyset$. 

The induced family $\mathfrak T_{M/\sbat} = \{ T_{S}/\sbat \tq S \in \stra_1 \}$ is a Thom-Mather system on the orbit space $M/\sbat$.
 The element $j\in \stres$ induces the map $j \colon M/\sbat \to M/\sbat$ preserving $\mathfrak T_{M/\sbat}$. So, it induces the map $j^* \colon \lau \Om * {\per p} {M/\sbat } \to \lau \Om * {\per p} {M/\sbat }$. Since $j^2 = -1\in \sbat$ then $j^*\circ j^*$ is the identity. So, the group $\zdos$ acts on $\lau \Om * {\per p} {M/\sbat }$ by $g \cdot \omega = j^* \omega$. 
 It also acts on  $\lau H * {\per p} {M/\sbat }$ by $g \cdot [\omega] = [j^* \omega]$. 
We now define
 $$
 \lau \Om * {\per p} {M/\sbat }^{-\zdos} = \set{ \omega \in \lau \Om * {\per p} {M/\sbat } | g\cdot \omega = -\omega}
 \quad\text{and}\quad
 \lau H * {\per p} {M/\sbat }^{-\zdos}= \set{ [\omega] \in \lau H * {\per p} {M/\sbat } | [g\cdot \omega] = -[\omega]}.
 $$
 As $\zdos$ is finite, the cohomology of $\lau \Om * {\per p} {M/\sbat }^{-\zdos}$ is easily seen to be
 $\lau H * {\per p} {M/\sbat }^{-\zdos}$.

\section{Invariant differential forms}
A key ingredient in this paper is the complex of $S^3$-invariant forms of $\lau \Om * {\per p} {M}$. It is a simpler sub-complex computing the same cohomology.

For the rest of this section we assume that the action  $\Phi \colon S^3 \times M \to M$ is a mobile action. In particular, the action of $S^3$ on $M\menos \Sigma$ is {\em almost free}, that is, the isotropy subgroup of any point of $M\menos \Sigma$ is finite.

\subsection{ The Lie algebra $\sudos$} \label{sudos}
Recall that identifying $\stres$ with the matrix group $\mathrm{SU}(2)$, its Lie algebra is $\sudos$. We shall consider $\{u_1,u_2,u_3\}$
an orthogonal  basis of $\sudos$, relatively to a bi-invariant metric $\kappa$ of $\stres$,  where
\begin{itemize}

\item $u_1$ generates the Lie algebra of the subgroup $\sbat$.

\item $[u_1,u_2] = u_3$, $[u_2,u_3] = u_{1}$,
$[u_3,u_1] = u_{2}$ and

\item  $\Ad(j)  \ u_1 = -u_1, \Ad(j) \ u_2 =u_2, \Ad(j)  \ u_3 = -u_3$. 

\end{itemize}

Consider the action $\Psi \colon \stres  \times \stres \to \stres $, defined by  
$\Psi(g,k) = k\cdot 
g^{-1}$. Associated with this action and $u_i$, for $i=1,2,3$, we have  the  fundamental  vector fields   $Y_i$ on $\stres$. 
They are
left invariant vector fields verifying $j_* Y_1 = -Y_1$, $j_*Y_2 = Y_2$ and $j_*Y_3 = -Y_3$.

 We shall write $\gamma_i \in \coho \Om 1 \stres$ the dual form of $Y_i$ relatively to the metric $ \kappa$: $\gamma_i =i_{Y_i} \kappa$, $i=1,2,3$. They are left invariant differential forms verifying
$
j^* \gamma_1 =- \gamma_1, j^*\gamma_2 = \gamma_2 \hbox{ and } j^*\gamma_3 = -\gamma_3.
$
 The differentials verify $d\gamma_1= \gamma_2 \wedge \gamma_3$, $d\gamma_2= -\gamma_1 \wedge \gamma_3$ and $d\gamma_3= \gamma_1 \wedge \gamma_2$.

\subsection{Fundamental vector fields  and characteristic forms.} Associated with the action $\Phi \colon \stres \times M \to M$ and the vector $u \in \sudos$, we have the {\em fundamental vector field } $X_u$ on $M$. 
For the sake of simplicity, we shall write
 $X_i = X_{u_i}$ with $i=1,2,3$. 
 This vector field  is defined on $M$ but we are going to work with its restriction to $M\menos \Sigma$. Since the action is mobile then these vector fields  are non-vanishing on $
M\menos \Sigma$. Moreover, the family $\{X_1(x),X_2(x),X_2(x)\}$ is a basis of the tangent space of the orbit $S^3(x)$ for any $x \in M\menos \Sigma$.
We have the equalities:
$
j_*X_{1} =j_* X_{u_1} = X_{\Ad(j)u_1} = - X_{u_1} = -X_1$
and 
$j_*X_{2} =X_2$, 
$j_*X_{3} =-X_3$ in the same way.

An {\em adapted metric} on $M\menos \Sigma$ is an $S^3$-invariant Riemannian metric $\mu$ on $M\menos \Sigma$ verifying
 \be\label{nu}
\mu (X_{v_1}(x), X_{v_2}(x)) =
\ \kappa(v_1,v_2) \hspace{2cm} \forall x \in M\menos \Sigma \hbox{ and } v_1,v_2 \in \sudos.
\ee
Such a metric always exists since the fundamental vector fields  are non-vanishing.

We denote by $\chii_{u} = i_{X_u} \mu\in \lau{\Om}{1}{}{ M\menos \Sigma}$ the
{\em characteristic form} associated
to $u \in \sudos$. Notice that, for each $g \in S^3$, we have
\be\label{coho}
g^*\chii_u = \chii_{\Ad(g^{-1})\cdot u}.
\ee
For the sake of simplicity, we shall write $\chii_i = \chii_{u_i}$, for $i=1 , 2,3$.
Since
$
L_{X_{u}}\chii_{v} = \chii_{[u,v]},
$
for each $u,v \in \sudos$, then we have
\be\label{lll}
\begin{array}{ll}
 L_{X_{1}}\chii_{1} = L_{X_{2}}\chii_{2}  = L_{X_{3}}\chii_{3}  
    = 0,  \hspace{2cm}& 
   L_{X_{1}}\chii_{2} = - L_{X_{2}}\chii_{1}  = - \chii_{3}\\[,4cm]
 L_{X_{1}}\chii_{3} = - L_{X_{3}}\chii_{1}  =  \chii_{2}& 
    L_{X_{2}}\chii_{3} = - L_{X_{3}}\chii_{2}  = - \chii_{1}.
    \end{array}
    \ee

Since
$
\chii_{k} (X_\ell) = \mu (X_\ell,X_k)= \delta_{\ell k},
$
each differential form $\om \in
\lau{\Om}{*}{}{ M\menos \Sigma}$ possesses a unique writing,
\begin{equation}\label{hamalau}
\om =\omega_0 +  
   \chii_{1} \wedge  \omega_1 +    \chii_{2} \wedge  \omega_2 +    \chii_{3} \wedge  \omega_3 +
 \chii_{1} \wedge  \chii_2 \wedge  \omega_{12} +     \chii_{1} \wedge  \chii_3 \wedge \omega_{13}+ 
  \chii_{2} \wedge  \chii_3 \wedge \omega_{23} + 
 \chii_1 \wedge \chii_2 \wedge \chii_3 \wedge \omega_{123},
\end{equation}
where the coefficients $ \ib  \om \bullet \in \coho \Om*{M\menos \Sigma} $ are  {\em horizontal forms}, that is, they verify  $i_{X_\ell}  \ib  \om \bullet =0$ for each $\ell=1, 2,3$. This is the {\em canonical decomposition} of $\om$.

The canonical decomposition of the differential of a characteristic form is
\be\label{ddd}
d\chii_{1} = e_1 + \chii_{2}  \wedge\chii_{3}  \hspace{2cm}
d\chii_{2} = e_2 -\chii_{1}  \wedge\chii_{3}  \hspace{2cm}
d\chii_{3} = e_3 + \chii_{1}  \wedge\chii_{2} 
    \ee
for some horizontal forms $e_1,e_2,e_3 \in \coho \Om 2 {M\menos\Sigma}$, called the {\em Euler forms}.
Notice that
\be\label{je}
j^*e_{1} = -e_1 \hspace{2cm}
j^*e_{2} =e_2 \hspace{2cm}
j^*e_{3} =-e_3 .
\ee

\subsection{Invariant differential forms}

A differential form $\om$ of $M\menos\Sigma$ is an {\em invariant form} when $g^*\omega =\omega$ for each $g \in\stres$ or, equivalently, $ \ib {L}{X_\ell}\om  =0$
for each  $\ell =1,2,3$. In fact, invariant differential forms are characterized by the following conditions:
\be\label{invll}
\begin{array}{lll}  \omega_0 \hbox{ and } \omega_{123} \hbox{ are basic forms} , 
&
L_{X_\ell}\omega_\ell=0, \ell =1,2,3 
 &
  L_{X_1} \omega_{23} =  L_{X_2} \omega_{13}= L_{X_3} \omega_{12} = 0 
   \\[,4cm]
 L_{X_1}\omega_2 = -L_{X_2}\omega_1 = -\omega_3 
 &
 L_{X_1}\omega_3 = -L_{X_3}\omega_1 = \omega_2,  
&
L_{X_2}\omega_3 = -L_{X_3}\omega_2 = -\omega_1,
 \\[,4cm]
L_{X_1} \omega_{13}  = L_{X_2} \omega_{23} = \omega_{12}
&
L_{X_1} \omega_{12}  = -L_{X_3} \omega_{23} = -\omega_{13}
 &
L_{X_2} \omega_{12}  = L_{X_3} \omega_{13} = -\omega_{23} 
\end{array}
\ee
(see \eqref{lll}).

The complex of invariant forms is denoted by $\coho \BOm* {M\menos \Sigma}$. The complex of invariant intersection differential forms is 
$\lau{\BOm}{*}{\per{p}}{M} = \coho \Om* {M\menos \Sigma} \cap \lau{\BOm}{*}{\per{p}}{M} $.
 Some cohomological computations are simplified by replacing the complex
 $\lau{\Om}{*}{\per{p}}{M}$
 by its subcomplex
 $\lau{\BOm}{*}{\per{p}}{M} $, since proceeding as in  \cite[Theorem I, p.~151]{MR0336651}, we have
 
\bp\label{inv}
The inclusion $\lau{\BOm}{*}{\per{p}}{M} \hookrightarrow \lau{\Om}{*}{\per{p}}{M} $ is a quasi-isomorphism for any perversity $\per p$.
\ep
Notice that $\lau \Om * {\per p} {M/\stres} = \{\omega \in \lau \BOm * {\per p} M \tq i_{X_\ell} \omega =0 \hbox{ for each } \ell = 1,2,3\}.
$

\subsection{Perverse degree of characteristic forms}
Notice that, for any singular stratum $S \in \strasing$, we have:
\be\label{99}
||\om||_S
=
\max \left\{ ||\omega_0||_S, 
   ||\chii_{\ell} \wedge  \omega_i ||_S,
 ||\chii_{\ell} \wedge  \chii_k \wedge  \omega_{\ell k} ||_S, 
|| \chii_1 \wedge \chii_2 \wedge \chii_3 \wedge \omega_{123}|| \tq 1 \leq \ell < k \leq 3
 \right\}
 \ee
 for any $\omega \in \coho \Om*{M\menos \Sigma}$ (cf. \eqref{PropGradPer}).
When $S \in \mathscr S_3$ is a fixed stratum then the orbits of the action are tangent to the fibers of $\tau_S \colon D_S \to S$. So, we have
\be\label{100}
 || \chii_{\ell} \wedge \alpha||_S -1
 =
 || \chii_{\ell} \wedge  \chii_{k}\wedge \alpha||_S-2
 =
 || \chii_{1} \wedge  \chii_{2}\wedge  \chii_{3}\wedge \alpha||_S-3
 =
 ||\alpha||_S,
 \ee
 for each $\alpha \in \coho\Om*{M\menos \Sigma}$ and each $1 \leq \ell < k \leq 3$.

  \smallskip

In order to control the perverse degree of characteristic forms  relatively to mobile strata  we need to enhance adapted metrics with richer properties.

\begin{definition}\label{def:adapted}
An adapted metric $\mu$ on $M\menos \Sigma$ is an {\em adjusted metric} if 
\be\label{battt}
\mu (X_v(x),w) =0
\ee
whenever
\begin{itemize}
\item  $x \in  \ib  D S \menos \Sigma$ for some $S \in  \stra_1$,
\item $w$ is a vector tangent to the fibers of  $\ib    \tau S \colon \ib  D S \menos \Sigma \to S$ at $x$, and 
\item $v \in \sudos$ belongs to the $\kappa$-orthogonal of  $\sudos_y$, Lie algebra of $S^3_{y}$ with $y=\ib  \tau S (x)$.
\end{itemize}

\end{definition}

\begin{proposition}
\label{batt}
Every mobile action admits an adjusted metric.
\end{proposition}
\begin{proof}
A convex combination of adapted metrics is an adapted metric. So, by using partitions of unity, we can reduce the problem to the following two cases:

\smallskip

$\bullet$ $M=T_S$ for some $S \in \stra_1$. In this case, $\Sigma=S$. We set $\mu'$ an adapted metric on $T_S\menos S$.

We put $\mathcal{K}$ (resp. $\mathcal G$) the sub-bundle of $T_S\menos S$ tangent to the fibers of $\tau_S$ (resp. the orbits of the action).  The bundle $ \mathcal{G} \cap \mathcal{K}$ is of constant rank equal to one. In fact, we have 
\bee
 \mathcal{G}_x \cap \mathcal{K}_x = \{X_v(x) \tq v \in  \sudos_y\}
 \eee
   for each $x \in T_S\menos \Sigma$ with $y=\tau_S(x)$. Let us consider the $S^3$-invariant decomposition
\bee\label{+++}
T\left(T_S\menos \Sigma\right)  = 
\mathcal{D} \oplus \mathcal{K} \oplus  \left( \mathcal{G} + \mathcal{K}\right)^{\bot_{\mu'}},
\eee
where $\mathcal{D}  = ( \mathcal{G} \cap \mathcal{K} )^{\bot_{\mu'}}  \cap\mathcal{G} $.
Since $\mu' =\kappa$ on $\mathcal{G}$ (cf. \eqref{nu}) then we have  $\mathcal{D}_{x} = \{X_v(x) \tq v \in  \sudos_y^\bot\}$ for each $x \in T_S\menos \Sigma$ with $y = \tau_S(x)$.
We denote by $\mu'_1, \mu'_2$, and $\mu'_3$ the restrictions of $\mu$ to each term of the above decomposition.
 The Riemannian metric $\mu$  defined by:
\bee
\mu = \mu'_1 +  \mu'_2 + \mu'_3,
\eee
is an adapted metric. It also satisfies  \eqref{battt} since $w \in \mathcal K_x$ and $X_v(x) \in \mathcal D_x$.

\smallskip

$\bullet$ $M=T_Q$ for some $Q \in \stra_3$. The open subset $T_Q \menos \Sigma$ is $\stres$-equivariantly diffeomorphic to $(D_Q\menos \Sigma )\times (0,\infty)$. The action of $\stres$ on $D_Q\menos \Sigma$ has no fixed points. The previous step gives an adjusted metric $\mu$ on  $D_Q\menos \Sigma$.
As the tubular neighborhood $T_S$ of any stratum $S \in \stra_1$ is the product $T_{S \cap (D_Q\menos \Sigma )} \times (0,\infty)$, the metric $\mu + dr^2$  is an adjusted metric on $T_Q\menos \Sigma$.
\end{proof}

For a such metric we can compute  the terms appearing in formula \eqref{99}.

\bp\label{35}
Let us suppose that $M\menos \Sigma$ is endowed with an adjusted metric.
Given a stratum $S \in \stra_1$ and a horizontal form $\alpha \in \coho\Om*{M\menos \Sigma}$ we have
 \be\label{101}
 || \chii_{\ell} \wedge \alpha||_S
 =
 || \chii_{\ell} \wedge  \chii_{k}\wedge \alpha||_S
 =
 || \chii_{1} \wedge  \chii_{2}\wedge  \chii_{3}\wedge \alpha||_S
 =
 ||\alpha||_S+1,
 \ee
 for each $1 \leq \ell < k \leq 3$.
\ep
\begin{proof}
Without loss of  generality, we can suppose $M=T_S$ and $\Sigma =S$. We proceed in two steps. 

\smallskip

\emph{Step $\leq$.} 
 Following \eqref{PropGradPer} it suffices to prove that 
  $|| \chii_\ell \wedge  \chii_k||_S \leq 
 1$ and 
 $|| \chii_1 \wedge  \chii_2 \wedge  \chii_3||_S \leq 
 1$.
 We deal with the first inequality, the second one can be approached in the same way.
 If $|| \chii_\ell \wedge  \chii_k||_S =2$ then there exists $x \in T_S \menos S$ and $v,w\in \mathcal K_x$ with 
 $\chii_\ell \wedge  \chii_{k}(v,w)\ne 0$. 
Denote by $\lfloor u_\ell,u_k\rfloor$ the vector subspace generated by $u_\ell$ and $u_k$.
Since $\dim \lfloor u_\ell,u_k\rfloor = 2 =\dim \sudos_y^\bot$ and $\dim \sudos =3$, then there exist $v_1 \in \sudos_y^\bot$ and $v_2 \in \sudos$ with $\chii_{v_1} \wedge  \chii_{v_2 }(v,w)\ne 0$.  
 This is impossible since $\chii_{v_1}(v) = \chii_{v_1}(w)=0$ (cf. \eqref{battt}).

\smallskip

\emph{Step $\geq$}. Since the result is clear for  $\alpha=0$, let us suppose $||\alpha||_S = a > 0$. 
So, there exist $x \in T_S\menos \Sigma$ and $\{w_0, \ldots,w_{a-1}\} \subset \mathcal K_x$,  with $\alpha(w_0, \ldots,w_{a-1},-)\ne 0$. Here, $\mathcal K$ denotes the sub-bundle of $T(T_S\menos S)$ tangent to the fibers of $\tau_S$ as defined on the proof of the previous Proposition.
Since the perverse degree is $\stres$-invariant, we can suppose that $\stres_{y=\tau_S(x)} \supset \sbat$
(cf. \eqref{PropGradPer}).
This gives $\sudos_y = \lfloor u_1 \rfloor $ and therefore $X_1(x) \in \mathcal K_x$.

The adjoint map associated to the group $\stres$ is the covering $\stres \to \mathrm{SO}(3)$. Since there exists a rotation sending $u_1$ to $u_\ell$ then there exists $g \in \stres$ with $\Ad(g)(u_1) =u_\ell$. Since the perverse degree is $S^3$-invariant  then it suffices to prove 
$ || \chii_{1} \wedge   \alpha||_S
 \geq  a+1$,
 $ || \chii_{1} \wedge  \chii_{u}\wedge \alpha||_S
 \geq  a+1$
 and
 $
 || \chii_{1} \wedge  \chii_{u}\wedge  \chii_{v}\wedge \alpha||_S
 \geq
 a+1,
$
where $u,v\in \sudos$
(cf. \eqref{PropGradPer} and \eqref{coho}).
Without loss of generality we can suppose that $u,v \in \lfloor u_1\rfloor^\bot$.
The inequality  comes from: 
\begin{eqnarray*}
0 &\ne &  \alpha(w_0, \ldots , w_{a-1}, -) =(\chii_1\wedge \chii_u \wedge \chii_v \wedge \alpha )(X_1(x),X_u(x),X_v(x), w_0, \ldots , w_{a-1}, -) \\
&= &
(\chii_1\wedge \chii_u\wedge \alpha )(X_1(x),X_u(x),w_0, \ldots , w_{a-1}, -) 
=
 (\chii_1 \wedge \alpha )(X_1(x),w_0, \ldots , w_{a-1}, -),
\end{eqnarray*}
since $\{w_0, \ldots,w_{a-1},X_1(x)\} \subset \mathcal K_x$.

If $a=0$, we just have $\alpha(-)\ne 0$  and the same argument applies by just omitting the vectors $w_0,\dots,w_{a-1}$.
\end{proof}

\subsection{Circle actions}

In \cite{MR1116314} a Gysin sequence is obtained for any nontrivial smooth circle action by doing a similar study and using more restrictive perversities. For the convenience of the reader, in this section we obtain that Gysin sequence for general perversities in a shorter way, using the techniques and presentation to be applied later for the case of mobile $S^3$-actions. 

Fix a nontrivial smooth action $\Psi \colon \sbat \times M \to M$.
  Here, we just have one fundamental vector field  $X$ and one characteristic form $\chii$ relatively to an adapted metric (the notion of adjusted metric does not apply here). This form is $\sbat$-invariant and verifies $||\chi||_S =1$ on $\stra_1$ (the family of fixed strata). Its differential $e=d\chi$, the {\em Euler form}, belongs to $\lau \Om 2 {\per e} {M/\sbat}$ where the {\em Euler perversity} $\per e$ is defined by $\per e(S) =2$ on $\stra_1$. We also use the {\em characteristic perversity } $\per \chii$ defined by $\per \chii(S)=1$ on $\stra_1$.

As $X$ is tangent to the fibers of $\tau_S\colon D_S\to S$ for any $S\in\mathscr{S}_1$, we have that for any $\alpha, \beta\in\lau\Om *{\rm bas} {M\menos \Sigma}$,   $||\alpha + \chi\wedge\beta||_S=\max\{ ||\alpha||_S , ||\chi\wedge\beta||_S\}$ holds. Using this fact, for any perversity $\per p$, the complex  $\lau \BOm*{\per p} M $ is
\be\label{decompeseuno}
\Set{ \alpha + \chii \wedge \beta |
	 \alpha \in \lau \Om *{\rm bas} {M\menos \Sigma}, 
\beta\in \lau \Om *{\per p - \per \chi} {M/\sbat} \hbox{ with }
\left\{
\begin{array}{l}
||\alpha|| \leq \per p(S), \hbox{ and }\\
 ||d\alpha +  e   \wedge \beta||_S\leq \per p(S) 
 \end{array}
 \right.
 \forall S \in \stra_1
}.
\ee
The \emph{ integration operator} is the differential operator $\fint \colon \lau \BOm * {\per p} M \to \lau \Om 
* {\per p  - \per \chi} {M/\sbat}$ defined by $\fint \omega = i_X \omega$, that is, $\fint (\alpha + \chi \wedge \beta) = \beta$.
 Considering the kernel and the image of $\fint$, we have the short exact sequence
 $$
 0 \to  \lau K* {\per p} M = \lau \Om *{\per p} {M/\sbat} \to  \lau \BOm * {\per p} M  \to  \lau I * {\per p} M \to 0,
 $$
 which induces the following Gysin sequence (see \cite{MR1116314}).
 \begin{proposition}\label{Gysin1}
 For each perversity $\per 0 \leq \per p \leq \per t$ we have the long exact sequence 
 $$
 \cdots \to \lau H *{\per p} {M/\sbat} \to \coho H*{M} \to \lau H {*-1}{\per p - \per e} {M/\sbat} \to
 \lau H {*+1}{\per p } {M/\sbat} \to \cdots.
 $$
 \end{proposition}
 \begin{proof}
 For each $\beta \in \lau \Om * {\per p - \per e} {M/\sbat} $ we have $\chi \wedge \beta \in \lau \BOm * {\per p } {M}$. Since $\fint (\chi \wedge \beta) = \beta$ then it  suffices to prove that the inclusion $I \colon \lau \Om * {\per p - \per e} {M/\sbat} \hookrightarrow  \lau I * {\per p} M$ is a quasi-isomorphism. 
 Notice that
 $$
  \lau I * {\per p} M = \set{ \beta \in \lau \Om * {\per p -\per \chi} {M/\sbat} | \exists \alpha \in \lau \Om *{\rm bas} {M\menos \Sigma} \hbox{ with } ||\alpha|| \leq \per p(S), \hbox{ and } 
 ||d\alpha +  e   \wedge \beta||_S\leq \per p(S)\ 
 \forall S \in \stra_1}.
 $$
 We proceed in three steps.

$\bullet$ {\em Step 1: $\stra_1=\emptyset$}. The action is almost-free. In this case $\Sigma= \emptyset$ and therefore 
$\lau \Om * {\rm bas} {M\menos \Sigma}  = \lau\Om * {\per p } {M/\sbat}  \subset  \lau I * {\per p - \per e } M \subset  \lau \Om * {\per p - \per \chi}  {M/\sbat} =  \lau \Om * {\rm bas} {M\menos \Sigma}$. 
In other words, the map $I$ itself  is an isomorphism.

\smallskip

 $\bullet$ {\em Step 2:  $M = T_S$ for some $S \in \stra_1$}. Recall that $\tau_S \colon T_S \to S$ is an $\sbat$-invariant smooth bundle whose fiber is $\R^{2b+2}$ for some $b \geq 0$. In fact,  the group $\sbat$ acts trivially on $S$ and orthogonally on the fiber $\R^{2b+2}$ having the origin as the only fixed point. Notice that the action of $\sbat$ on the unit sphere $S^{2b+1}$ is almost-free.
 
 Consider a good covering $\mathcal U$ of $S$ and $\{f_U \tq U\in \mathcal U\}$ a subordinated partition of unity. The family $\{\tau_S^{-1}(U),\tq u \in \mathcal U\}$ is an open covering of $T_S$ having 
 $\{f_U \circ \tau_S  \tq U\in \mathcal U\}$ a subordinated partition of unity. These maps are $\sbat$-invariant smooth maps constant on the fibers of $\tau_S$. This last property implies that $||f_U \circ \tau_S ||_S =||d(f_U \circ \tau_S )||_S=0$. So, the covering $\mathcal U$ possesses a subordinated partition of unity living in $\lau\BOm*{\per 0} M$.
 
By applying Bredon's trick \cite[p.~289]{MR1700700}, we can reduce the problem to the case where $M = \R^{\dim S} \times \R^{2b+2}$, with $\tau_S$ being the projection onto the first factor. The action of the group $\sbat$ is trivial on the first factor.

 Contracting this factor to a point, we reduce the problem to the case $M= \R^{2b+2} = \rc S^{2b+1} = (S^{2b+1} \times [0,\infty))/(S^{2b+1}  \times \{0\})$, where $S^1$ acts on the factor $S^{2b+1}$.
 Here, the stratum
 $S $ is the apex of the cone. We have $\per \chi (S) = 1$ and $\per e(S)=2$. 
 The number $p \in \per \Z$ is defined by $\per p(S) =p$.
 We need to prove that the inclusion
 \be\label{ayuda}
 I \colon \lau \Om * {\per p - \per e} {\rc S^{2b+1}/\sbat} \hookrightarrow
\lau I * {\per p } {\rc S^{2b+1}} = \{ i_Z \omega \tq \omega \in \lau \BOm {*+1} {\per p} {\rc S^{2b+1}}\}
 \ee
 is a quasi-isomorphism. 
Notice first that
 $$
  \begin{array}{lcl}
 \lau \Om {*<p-2 } {\per p - \per e} {\rc S^{2b+1}/\sbat} &=& \lau \Om {*<p-2}{\rm bas} { S^{2b+1}  \times (0,\infty)}, \\[,2cm]

 \lau \Om {p -2}  {\per p - \per e} {\rc S^{2b+1}/\sbat} &=& \Set{ \beta \in  \lau \Om {p-2} {\rm bas} { S^{2b+1}  \times (0,\infty)} | d\beta \equiv 0 \hbox{ on }  S^{2b+1}  \times (0,2)}, \\[,2cm]
 
\lau \Om {*> p -2}  {\per p - \per e} {\rc S^{2b+1}/\sbat} &=&  \lau \Om {*>p-2} {\rm bas} { S^{2b+1} \times (0,\infty), S^{2b+1}  \times (0,2)}, \\[,2cm]

 \lau I {* <p-1}{\per p}{\rc S^{2b+1}}&=& \Set{  i_X\omega | \omega \in \lau  \BOm {* + 1 <p}{} {S^{2b+1}\times  (0,\infty)}} =_{(1)}    \lau  \Om {* <p-1}{\rm bas} {S^{2b+1}\times  (0,\infty)},  \\[,2cm]
 
 \lau I {p-1}{\per p }{\rc S^{2b+1}}  &=& \Set{  i_X \omega | \omega \in \lau  \BOm {p}{} {S^{2b+1}\times  (0,\infty)}  \text{ and } d\omega \equiv 0 \hbox{ on }  S^{2b+1}  \times (0,2)}, \hbox{ and} \\[,2cm]
 
 \lau I {*> p-1}{\per p }{\rc S^{2b+1}}  &=& \Set{  i_X \omega | \omega \in \lau  \BOm {*+1>p}{} {S^{2b+1}\times  (0,\infty),S^{2b+1}\times  (0,2)} }\\[,1cm] 
 &=_{(2)} & \lau \Om {*>p-1} {\rm bas} { S^{2b+1}  \times (0,\infty), S^{2b+1}  \times (0,2)} ,
 \end{array}
 $$
 where $=_{(1)}$ is given by the previous step and $=_{(2)}$ comes from the fact that $\beta \equiv 0$ on $ S^{2b+1} \times (0,2)$ implies $\omega = \chi \wedge \beta \equiv 0$ on $S^{2b+1} \times (0,2)$.
 Since $ \lau \Om {p -2}  {\per p - \per e} {\rc S^{2b+1}/\sbat}  \cap d^{-1}(0) =\lau I {p-2}{\per p }{\rc S^{2b+1}}  \cap d^{-1}(0) $ then it suffices to study the degrees $*\geq p-1$. 

\begin{itemize}
\item[] \fbox{$*=p-1$}
 Since $\lau H {p-1} {\per p - \per e} {\rc S^{2b+1}/\sbat} =0$ then we need to prove $\coho H {p-1}{ \lau I * {\per p} {\rc S^{2b+1}} }=0$, that is:
  $$
  \left\{
  \begin{array}{l}
  \omega \in \coho \BOm p   {S^{2b+1} \times (0,\infty)} \\[,2cm]
  \hbox{with }d\omega \equiv 0 \hbox{ on } S^{2b+1} \times (0,2) \hbox{ and } 
  d i_X \omega =0
  \end{array}
  \right.
  \Longrightarrow
    \left\{
  \begin{array}{l} \exists \eta \in \coho \BOm {p-1}  {S^{2b+1} \times (0,\infty)} \\[,2cm]
  \hbox{with }  i_X  \omega = d i_X  \eta.
    \end{array}
  \right.
  $$
If $p=0$ then we can consider $\eta=0$ since $\omega $ is constant. Let us  suppose $p\geq 1$. Consider $\eta' = \int_1^- \omega$. It is an element of $\coho \BOm {p-1}  {S^{2b+1} \times (0,\infty)} $ since the action of $\sbat $ on the $(0,\infty)$-factor is trivial.
   A straightforward calculation gives $\omega = \pr^* \omega(1) + d \eta'  + \int_1^- d\omega$. 
  Here, $\omega(1)$ is the restriction of $\omega$ to $S^{2b+1} \times \{1 \}$ and $\pr \colon S^{2b+1} \times (0,\infty) \to S^{2b+1} \times \{ 1\}$ is the map defined by $\pr(x,t) =(x,1)$.  A straightforward calculation gives
  $$
i_X \omega = i_X \pr^* \omega(1) + i_X d \eta'+ i_X\int_1^-d\om =  i_X \pr^*\omega(1) -  d i_X \eta' - \int_1^- di_X \om =  i_X\pr^*\omega(1) -  d i_X \eta' .
  $$ 
  By hypothesis the differential form $\pr^*\omega(1)$ is a cycle of $\lau \BOm p {} {S^{2b+1} \times (0,\infty)}$.  Condition $  \per p \leq \per t$ implies $p  \leq \codim S - 2 = 2b$. Since $p \geq 1$, $\pr^*\omega(1)$ is exact and so, there exists $\eta'' \in \coho \BOm {p-1} {S^{2b+1} \times (0,\infty)}$ with $\pr^*\omega(1) = d \eta''$. We end the proof taking $\eta = -\eta' - \eta''$.
  
\item[] \fbox{$* \geq p$}
 Since $\lau H {* \geq p} {\per p - \per e} {\rc S^{2b+1}/\sbat} =0$ then we need to prove $\coho H {* \geq p}{ \lau I * {\per p} {\rc S^{2b+1}} }=0$, that is:
  $$
  \left\{
  \begin{array}{l}
  \omega \in \coho \BOm {* +1\geq p+1}   {S^{2b+1} \times (0,\infty),S^{2b+1} \times (0,2)} \\[,2cm]
  \hbox{ with } d i_X \omega =0
  \end{array}
  \right.
  \Longrightarrow
    \left\{
  \begin{array}{l} \exists \eta \in \coho \BOm {* \geq p}  {S^{2b+1} \times (0,\infty)} \\[,2cm]
   \hbox{with }d\eta \equiv 0 \hbox{ on } S^{2b+1} \times (0,2)\hbox{and }  i_X  \omega = d i_X  \eta.
    \end{array}
  \right.
  $$
  The proof of the previous case applies verbatim with $\omega(1)=0$.
\end{itemize}
 
 $\bullet$ {\em Final Step}. Consider the invariant open covering $\mathcal V =\{T_S \tq S \in \stra_1\} \sqcup \{M\menos \Sigma\}$ of $M$. 
 We fix a smooth map $\lambda \colon [0,\infty) \to [0,1]$ verifying $\lambda=1$ on $[0,2]$ and $\lambda=0$ on $[3,\infty)$. 
 The map $f_S\colon T_S \to [0,\infty)$ is defined by $f_S(x) = \lambda(
 \rho_S(x))$ (see \subsecref{TMs}). It is an $\sbat$-invariant smooth map, constant on the fibers of $\tau_S \colon D_S \to S$, which gives  $||f_S||_S = ||df_S||_S=0$. So, the family $\{ f_S \tq S \in \stra_1\} \sqcup \{ 1 - \sum f_S\}$ is a partition of unity, subordinated to $\mathcal V$, living in $\lau\BOm*{\per 0} M$.
 Now, it suffices to apply   Bredon's trick \cite[p.~289]{MR1700700} and the previous cases.
 \end{proof}

\begin{remarque}  \label{rem}
In fact, we have proved that the operator 
 \bee\label{iso}
 f \colon ( \lau \Om * {\per p} {M/\sbat} \oplus \lau \Om {*-1} {\per p - \per e} {M/\sbat},D) \TO (\lau \BOm*{\per p} M,d),
 \eee
 defined by $f(\alpha,\beta) = \alpha + \chi \wedge \beta$, where $D(\alpha , \beta) = (d\alpha + e \wedge \beta, -d\beta)$, is a quasi-isomorphism. 
%
 \end{remarque}

\subsection{Twisted product} \label{twist}
The model  of the tubular neighborhood of a semi-mobile stratum is given by twisted products. We present this notion. First of all, we consider a smooth action $\Theta \colon N \times E \to E$
 of the normalizer $N$ on  a manifold $E$.
 It induces the action $\Phi \colon \stres \times \left(S^3 \times_N E \right)\to \left(S^3 \times_N E\right)$ defined by 
 $
 g \cdot  \left< k,x\right> =  \left<g \cdot k,x\right> .
 $
 
 This action is a mobile action and the strata are of the form $\stres \times_N S$ with $S \in \stra$, stratification induced by $\Theta$.
A perversity $\per p$ on $E$  determines a perversity on  the twisted product, still denoted $\per p$, defined by 
$\per p\left(\stres \times_N S\right) = \per p(S)$.

Given an $N$-invariant Thom-Mather system $\mathfrak T_E = \{ T_S \tq S \in \stra_1\}$ on $E$ we can consider the following Thom-Mather systems 
\begin{itemize}
\item
$\mathfrak T_{\stres\times E}  = \{ \stres \times T_S \tq S \in \stra_1\}$ on the product, and
\item
$\mathfrak T_{\stres\times_N E}  = \{ \stres \times_N T_S \tq S \in \stra_1\}$ on the twisted product.

\end{itemize}
Relatively to these Thom-Mather systems we have the equality
\begin{equation}\label{perdegig}
||\Pi^*\omega||_{\stres \times S} = ||\omega||_{\stres \times_N S}
\end{equation}
where $S \in \stra_1$, $\omega \in \lau \Om *{} {\stres\times_N (E\menos \Sigma)}$ and $\Pi \colon \stres \times E \to \stres \times_N E $ is the canonical projection. This map is an $N$-bundle and verifies 
$\pi(g,x) = \pi (g \cdot h^{-1},h\cdot x)$ where $(g,x) \in \stres \times E$ and $h \in N$.

The goal of this section is to write the $\stres$-invariant intersection forms of the twisted product $S^3 \times_N E $ in terms of the intersection forms of $E$. 

First we establish some notation.

\medskip

\begin{enumerate}[(a)]

\item $Y_u \in \mathfrak{X}(S^3)$ is the fundamental vector field  associated to $u\in \sudos$ relatively to the right action: $S^3 \times S^3 \to S^3; (g,k) \mapsto k \cdot g^{-1}$. It is a left invariant vector field . For the sake of simplicity we shall write $Y_{u_\ell}=Y_{\ell}$ for $\ell \in \{1,2,3\}$.

\item $Z \in \mathfrak{X}(E)$  is the fundamental vector field   of the action: $\Psi \colon  \sbat \times E \to E$, induced by $\Theta$. It verifies $j_*Z = -Z$. Let $\epsilon$ be an $N$-invariant metric on $E$. The \emph{characteristic  form} $\zeta = \iota_Z\epsilon $ verifies $j^* \zeta=-\zeta$ and the associated \emph{Euler form} $e = d\zeta$ verifies $j^*e = e$.

\item Let  $\gamma_u \in \coho {\Om} 1 \stres$ be the 
dual form of $Y_u$, that is, $\gamma_u = i_{Y_u} \kappa$, $u\in\sudos$.  Notice that $\kappa(u,v) = 
\gamma_{u}(Y_v) $.
These forms are invariant by the left
action of $\stres$. For the sake of simplicity we shall write $\gamma_{u_k} = \gamma_k$ for $\ell \in \{1,2,3\}$.
They verify $L_{Y_1} \gamma_1=0,  L_{Y_1} \gamma_2=-\gamma_3$, $L_{Y_1} \gamma_3=\gamma_2$, $d\gamma_1=\gamma_2\wedge \gamma_3$, $d\gamma_2 = -\gamma_1 \wedge \gamma_3$ and $d\gamma_3 - \gamma_1 \wedge \gamma_2$ 
(cf. \eqref{lll}). 

\item The group $N$ acts on the complex of differential forms $\coho \Om * \stres$ by the left. So, the group $\zdos = N /\sbat$ acts on the complex of $\sbat$-left invariant forms of $\stres$, which is $\bigwedge^* (\gamma_1,\gamma_2,\gamma_3)$.  This action is given by 
\be\label{j}
g \cdot \gamma_\ell = (-1)^\ell \gamma_\ell,
\ee
for  $\ell=1,2,3$.

 \end{enumerate}

\bp
\label{twist2} The canonical projection $\Pi \colon \stres \times E \to \stres \times_N E $ induces the identification
\begin{equation}\label{identi}
\lau \BOm * {\per{p}} {S^3 \times_N E}
=
\Set{ \om \in  \displaystyle {\bigwedge}^*( \gamma_{1}, \gamma_2, \gamma_{3}) \otimes \lau \Om *{\per{p}} E |
\begin{array}{l}
i_{Y_1}\om = - i_{Z}\om \\
L_{Y_1}\om =-  L_{Z}\om
\end{array}} ^{\zdos}
\end{equation}

\ep
\begin{proof}
Since the map $\Pi$ is $\stres$-invariant then $\Pi^*$ induces a monomorphism between 
$$
 \coho\BOm * {\stres \times E} = \{ \omega \in \coho\Om * {\stres \times E} \tq g^*\omega = \omega \ \forall g \in \stres\}  = {\bigwedge}^* (\gamma_1,\gamma_2,\gamma_3) \otimes  \coho \Om* E  $$
and
$$  \coho\BOm * {\stres \times_N E} = \{ \omega \in \coho\Om * {\stres \times_N E} \tq g^*\omega = \omega \ \forall g \in \stres\} .
$$
So, we can identify $ \coho\BOm * {\stres_N \times E}$ with 
 $$
\Set{ \om \in  {\bigwedge}^*( \gamma_1, \gamma_2, \gamma_3) \otimes \lau \Om *{} {E} | 
\begin{array}{l}
i_{Y_1}\om =  - i_{Z}\om \\
L_{Y_1}\om = - L_{Z}\om \\
g^*\om = \om 
\end{array}
} 
=
\Set{ \om \in  {\bigwedge}^*( \gamma_1, \gamma_2, \gamma_3) \otimes \lau \Om *{} {E} |
\begin{array}{l}
i_{Y_1}\om =  - i_{Z}\om \\
L_{Y_1}\om = - L_{Z}\om  
\end{array}
}^{\zdos}.
 $$
A similar identification is obtained for $E\menos \Sigma_E$ instead of $E$. Using \eqref{perdegig} we get  \eqref{identi}.\end{proof}

\begin{corollaire}\label{simpli}
We have the identification 
\bee\label{forsimpli}
\lau \BOm * {\per{p}} {S^3 \times_N E} =
 \lau \BOm* {\per{p}}{E }^{\zdos}
 \oplus \lau \BOm{*-2} {\per{p}}{E}^{-\zdos}
 \oplus
\{ \xi \in  
\lau \Om{*-1} {\per{p}}{E}^{-\zdos} \tq
\ib L {Z} \ib L {Z} \xi=-\xi\},
\eee
where the differential becomes $D'(\alpha, \beta, \xi)= (d\alpha, d\beta - i_Z\alpha,-d \xi)$. The third term of this direct sum is acyclic.
\end{corollaire}
\begin{proof}
\propref{identi} shows that an element of $\lau \BOm * {\per{p}} {S^3 \times_N E}$ is of the form
 $$
 \alpha   -  \gamma_1 \wedge i_Z\alpha   + \gamma_2 \wedge L_{Z} \xi + \gamma_3 \wedge \xi 
   +
 \gamma_1 \wedge \gamma_2 \wedge i_{Z} L_{Z} \xi 
   + \gamma_1 \wedge \gamma_3 \wedge  i_{Z} \xi
   +
   \gamma_2 \wedge \gamma_3 \wedge \beta  
   - \gamma_1 \wedge \gamma_2 \wedge \gamma_3 \wedge i_{Z}\beta,
 $$
with $\alpha \in \lau \BOm* {\per{p}}{E}$,  $ \beta \in \lau \BOm{*-2} {\per{p}}{E}$, $\xi \in \lau \Om{*-1} {\per{p}}{E}$, $g \cdot \alpha = \alpha$, $g \cdot \xi = -\xi$, $g \cdot\beta = -\beta$, and $\ib L {Z} \ib L {Z} \xi=-\xi$.
For the calculation of $D$ we compute the differential of the previous expression:
$$
 d\alpha   -  \gamma_1 \wedge i_Z d\alpha   - \gamma_2 \wedge L_{Z} d\xi - \gamma_3 \wedge d \xi 
   -
 \gamma_1 \wedge \gamma_2 \wedge (i_{Z} L_{Z}d\xi +\xi) 
   + \gamma_1 \wedge \gamma_3 \wedge (L_Z\xi - i_{Z} d\xi)
   +
   $$
   $$
   \gamma_2 \wedge \gamma_3 \wedge(d \beta  -i_Z\alpha) 
   + \gamma_1 \wedge \gamma_2 \wedge \gamma_3 \wedge i_{Z}d\beta.
 $$

We verify the acyclicity property.
Let  $\xi\in \lau \Om * {\per{p}} {E}^{-\zdos}$ be a cycle. The differential form  $\eta = i_ZL_Z \xi \in \coho \Om *  {E \menos \Sigma} $ verifies
\begin{itemize}
\item $g \cdot \eta = j^*\eta = j^* i_ZL_Z \xi = -i_Zj^*L_Z \xi = i_ZL_Z j^*\xi = -i_ZL_Z \xi = -\eta$.
\item $L_ZL_Z\eta =L_ZL_Z i_ZL_Z \xi = i_ZL_ZL_Z L_Z \xi = -i_ZL_Z \xi = -\eta$.


\item  $d\eta = di_ZL_Z \xi = L_ZL_Z \xi = -\xi$, since $d \xi =0$.
\item If $Q$ is a singular stratum of  $ \ib E S$ we have
$
||\eta||_{Q} = ||i_ZL_Z \xi||_{Q} \stackrel{\eqref{PropGradPer}}{\leq} ||\xi||_{Q}\leq  \per{p}(Q).
$
\end{itemize}
So, the complex $\{ \xi \in  
\lau \Om{*-1} {\per{p}}{E}^{-\zdos} \tq
\ib L {Z} \ib L {Z} \xi=-\xi\}$ is acyclic.
 \end{proof}
The following calculations  will be used in the next section. We use the $N$-action presented in \secref{circuloN}.

 \begin{corollaire}\label{lema}
 Let $\Theta \colon N \times M \to M$ be a smooth action. We have
 $$
 \lau H*{\per p /( \per p -\per e)} {M/\sbat}^{-\zdos} 
 =
\coho H {*-2q} {\Sigma}^{-(-1)^{q }\zdos}.
 $$
 where  
 \begin{itemize}
  \item   $\per p$ is a constant perversity $p$ on $M/\sbat$ verifying $\per 0 \leq \per p \leq \per t$, 
 \item $q$ denotes the integer part of $p/2$, and
  \item $\Sigma$ is the singular part of the induced $\sbat$-action, 
\end{itemize}
 \end{corollaire}
 \begin{proof}
 Consider an $N$-invariant Thom-Mather system (cf. \secref{circuloN}).
Using Mayer-Vietoris, we can suppose $M= T_S \cap T_{j(S)}$, for some $S \in \stra_1$. We write $\tau = \tau_S \cup \tau_{j(S)}$.
Since $\Sigma=S \cup j(S)$ is $\zdos$-invariant, then we can consider the operator 
\bee\label{JJ}
J \colon   \coho \Om  {*-2q } {\Sigma}^{-(-1)^{q}\zdos} \longrightarrow  \lau \Om* {\per{p}/(\per{p}-\per{e})} {M/\sbat}^{-\zdos} ,
\eee
 defined by $J (\alpha) = <\tau_S ^*\alpha \wedge e^{q}>$, where $e$ is the Euler form of the induced $\sbat$-action relatively to an $N$-invariant Riemannian metric on $M$. This metric always exists since the group $N$ is compact.
It is a well defined differential operator since 
\begin{itemize}
\item $
|| \tau ^*\alpha \wedge e^{q}||_{S} \leq  || \tau^* \alpha ||_{S}  + || e^{q}||_{S} \leq 
0+|| e^{q}||_{S}
 \leq 
2{q} 
\leq \per p(S),
$
and similarly $|| \tau ^*\alpha \wedge e^{q}||_{S} \leq \per p(j(S))$,

\item $d( \tau^*d\alpha \wedge e^{q} )= \tau^*d\alpha \wedge e^{q} $ and the operator is a differential operator,
\item 
$ i_X (\tau^*\alpha \wedge e^{q} ) =  i_X d(\tau^*\alpha \wedge e^{q} )=0$, which gives $\tau^*\alpha \wedge e^{q}
\in \lau \Om *{\rm bas} {M\menos S}$, and
\item  $j^* e =-e$.
\end{itemize}
This last property comes from  $j_*X =-X$, where $X$ is the fundamental vector field  of the circle action $\Psi$.

We claim that the operator $J$ is a quasi-isomorphism. Proceeding as in the proof of the \propref{Gysin1} we can reduce the question to the case where the stratum $S$ (resp. $j(S)$) is the apex $\tv$ (resp. $\tw$) of the cone $T_S= \rc_\tv S^{2b+1}$ (resp. $T_{j(S)}= \rc_{\tw} S^{2b+1}$) where the circle $\sbat$ acts orthogonally and almost freely on the sphere $S^{2b+1}$. We also have $\lau H*{\rm bas}{S^{2b+1}} = \coho H * {\C P^{b}}$.
Recall that $\per e (S)=2$ and $\per p(S)=p$. We distinguish two cases:

\begin{itemize}
 \item[] \fbox{$* \leq p-2$ or $*> p$ }.  We have the equality $\lau \Om {*} {\per p - \per e} {\rc_\tv S^{2b+1}/\sbat}  = \lau \Om {*} {\per p } {\rc_\tv S^{2b+1}/\sbat} $ and the analogous one for the apex $\tw$.  Therefore, and for degree reasons, we get 
 $$
 \lau \Om {*} {\per p/(p - \per e)} {(T_S \cup T_{j(S)})/\sbat}= 0 = 
 \coho \Om {*-2q}{\{ \tv, \tw\}}^{-(-1)^{q}\zdos}=  \coho \Om {*-2q}{S\cup j(S)}^{-(-1)^{q}\zdos}.
 $$
 
\item[] \fbox{$* = p -1$  or $p$}. We first assume $*=2q$. On one hand, we have
$$
   \coho H  {*-2q } {  S \cup j(S) }^{-(-1)^{q}\zdos} =
\left\{
\begin{array}{ll}
\coho H  {0} {  S }^{-(-1)^{q}\zdos}  = 0 \hbox{ (if $q$ even) or }  \R \hbox{ (if $q$ odd)} & \hbox{if } j(S) =S\\
\coho H  {0} {  S\cup j(S) }^{-(-1)^q\zdos} =\langle(1,(-1)^{q+1})\rangle=  \R    & \hbox{if } j(S)  \cap S =\emptyset.
\end{array}
\right.
$$
On the other hand, the typical conical calculations of the intersection cohomology (see, for example, \cite[Proposition 3.1.1]{MR2210257})  and the sequence \eqref{pasito} give $\lau H {*} {{\per p}/({\per p} - {\per e})} {\rc_\tv S^{2b+1}/\sbat}
= \coho H{2{q}}{\C P^{b}}$, which is generated by $e^{q}
$, since $0 \leq p = \per p(S) \leq \per t (S) =  2b$. The same formula is obtained for $\tw$. We have
$$
   \lau H* {\per{p}/(\per{p}-\per{e})} {(\rc_\tv S^{2b+1} \cup\rc_\tv S^{2b+1}) /\sbat}^{-\zdos} = 
\left\{
\begin{array}{ll}
 \coho H{2{q}}{\C P^{b}} ^{-\zdos} = 0 \hbox{ (if $q$ even) or }  \R \hbox{ (if $q$ odd)} & \hbox{if  } j(S) =S\\
\coho H{2{q}}{\C P^{b}\cup \C P^{b}}^{-(-1)^q\zdos} = \R & \hbox{if  } j(S)  \cap S =\emptyset,
\end{array}
\right. 
$$
where, in the case where $S\cap j(S)=\emptyset$, a generator of the cohomology is the class $([e]^q,(-1)^{q+1}[e]^q)$.

Finally, if $*$ is odd, both cohomology groups to be checked equal vanish because they involve classes in odd degrees of points or $\mathbb{C}P^b$.
\end{itemize}
 We get that $J$ is a quasi-isomorphism.
\end{proof}

\section{The integration operator $\fint$.}\label{OpInt}

The main tool we use in this work is the integration operator 
$$
\fint \colon \lau{\BOm}{*}{\per{p}}{M} \TO \lau{\Om}{*-3}{\per{p} - \per \chi}{M/\stres },
$$
defined by 
$\fint  \om   = (-1)^{\deg \om}  \  \ib  {i}{X_3}  \ib  {i}{X_2}  \ib  {i}{X_1} \om$, where $\per \chi$ is the  {\em characteristic perversity} defined by $
 \per \chi(S)= \left\{
 \begin{array}{ll}
 1 & \hbox{if  } S \in \stra_1 \\
 3 & \hbox{if  } S \in \stra_3 
 \end{array}
 \right.
 $. 
The operator $\fint$ is  a well defined differential operator since
\begin{itemize}
\item
$L_{A} i_{B} = i_{B} L_{A} +
 i_{[A,B]}$ when $A,B$ are vector fields  of $M\menos \Sigma$  and 
 \item $ \per{p}(S) \geq ||\om||_S \stackrel{\eqref{hamalau},\eqref{99}}{\geq} ||\chii_1 \wedge \chii_2 \wedge \chii_3 \wedge \ib  {i}{X_3}  \ib  {i}{X_2}  \ib  {i}{X_1}\om||_S \stackrel{\eqref{100}, \eqref{101}}{=\!=} ||  \ib i{X_3}  \ib  {i}{X_2}  \ib  {i}{X_1}\om||_S + \per \chii(S)= || \fint\om||_S + \per \chii(S)$ for each $S \in \stra$,
 \end{itemize}
 where we have  considered an adjusted metric $\mu$ on $M\menos \Sigma$. 
 
 The goal of this section is the computation of the cohomology of the complexes  $ \lau I * {\per p} M$ and $ \lau K* {\per p} M$.
 For the first one, we need to introduce the {\em Euler perversity} $\per e$, defined by
 $
 \per e = \per \chi +\per 1.
 $
 
 For the sake of simplicity we shall write
 \begin{eqnarray*}
  \ker \fint& =& \lau K * {\per p} M = \left\{ \omega \in \lau \BOm * {\per p} M \tq i_{X_3} i_{X_2}i_{X_1} \omega=0 \right\}, \text{ and}\\
 \Ima \fint & =& \lau I * {\per p} M = \left\{ i_{X_3} i_{X_2}i_{X_1} \omega  \tq \omega \in \lau \BOm {*+3} {\per p} M\right\}.
 \end{eqnarray*}
 \begin{proposition}\label{31}
 Let $\per p \leq \per t$ be a perversity on $M$. The natural inclusion 
 $\displaystyle
 I \colon \lau \Om * {\per p -\per e} {M/\stres} \hookrightarrow  \lau I * {\per p} M
 $
 is a quasi-isomorphism.
 \end{proposition}
 \begin{proof} 
 The inclusion makes sense if we prove that $\chii_1 \wedge  \chii_2 \wedge \chii_3 \wedge \alpha \in \lau \BOm * {\per p} M$ for each $\alpha \in \lau \Om {*-3} {\per p - \per e} {M/\stres}$. This comes from
 \begin{eqnarray*}
 L_{X_\ell}(\chii_1 \wedge  \chii_2 \wedge \chii_3 \wedge \alpha) &=&0 \hbox{ for each } \ell \in \{1,2,3\}
 \hbox{  (cf. \eqref{lll})}, \\
 \hbox{and for each $S \in\stra_1 \sqcup \stra_3$:}&&\\
 ||\chii_1 \wedge \chii_2 \wedge \chii_3 \wedge \alpha||_S &\stackrel{\eqref{100},\eqref{101}}{=}&
  ||\alpha||_S + \per \chii (S) \leq \per{p}(S)-\per{e}(S) +\per \chii(S) \leq \per{p}(S) \\
||d(\chii_1 \wedge \chii_2 \wedge \chii_3 \wedge \alpha)||_S &\stackrel{\eqref{hamalau},\eqref{99}}{\leq}  & \max(||\chii_1 \wedge \chii_2 \wedge \chii_3 \wedge d \alpha||_S ,|| d(\chii_1 \wedge \chii_2 \wedge \chii_3) \wedge \alpha||_S) \leq \max( \per{p}(S), ||\alpha||_S\\
&&+||d(\chii_1 \wedge \chii_2 \wedge \chii_3)||_S )\stackrel {\eqref{lll}} \leq 
\max( \per{p}(S), \per{p}(S) - \per{e}(S) + ||\chii_1 \wedge \chii_2 \wedge \chii_3||_S +1)\\
&&
\stackrel {\eqref{100},\eqref{101}} \leq \max( \per{p}(S), \per{p}(S) - \per{e}(S) + \per e(S))
= \per p(S)
\end{eqnarray*}
(cf.  \eqref{PropGradPer}).

In order to prove that $I$ is a quasi-isomorphism, we proceed in several steps. We use the $\stres$-invariant Thom-Mather system of \secref{TMs}.

$\bullet$ {\em Step 1: $\stra_1=\stra_3=\emptyset$}.  The action is almost-free. In this case $\Sigma= \emptyset$ and therefore 
$\lau \Om * {\rm bas} {M}  = \lau\Om * {\per p } {M/\stres}  \subset  \lau I * {\per p - \per e } M \subset  \lau \Om * {\per p - \per \chi}  {M/\stres} =  \lau \Om * {\rm bas} {M}  $. In other words, the map $I$ itself  is an isomorphism.

\smallskip

 $\bullet$ {\em Step 2:  $M = T_S$ for some $S \in \stra_1$}.
We have seen that $S = \stres \times_N S^{\sbat}$. The restriction $\tau_S \colon E_S = \tau_S^{-1}(S^{\sbat})\to S^{\sbat}$ is an $N$-invariant bundle. Notice that $E_S$ is a filtered space whose singular strata are the connected components of $S^{\sbat}$.
The fiber of this bundle is $\R^{2b+2}$, for some $b\in \N$. 
The group $\sbat$ acts trivially on $S^{\sbat}$  and   orthogonally on $\R^{2b+2}$ having the origin as the only fixed point.
Notice that the action of $\sbat$  on the unit sphere $S^{2b+1}$ is almost-free.

We identify $T_S$ with the twisted product $\stres \times_N E_S$ and we use the calculations of \secref{twist}. 
We have
 \be\label{idde2}
 \lau \Om * {\per p} {M/\stres} = \lau \Om * {\per p} {E_S/N} = \lau \Om * {\per p} {E_S/\sbat}^{\zdos}
 \ee
 for any perversity $\per p$.
 The integration $\fint$ becomes the map
\be\label{intedef}
\fint \colon \lau \BOm * {\per{p}} {S^3 \times_N E_S} \TO  \lau \Om{*-3} {\per{p} - \per \chi}{E_S/\sbat}^{\zdos} ,
\ee
defined by $(\alpha,\beta,\xi) \mapsto - i_Z \beta$ 
 (cf. \corref{simpli}). 
 The map $I$ becomes the inclusion 
 \be\label{inte}
 I \colon
 \lau \Om {*} {\per p -\per e} {M/\stres} =
 \lau \Om {*} {\per p -\per e} {E_S/\sbat}^{\zdos}  \TO  \lau I * {\per p} M   =
  \left\{ i_Z \beta \tq  \beta \in \lau \BOm{*+1} {\per{p}}{E_S}^{-\zdos} \right\}
 \ee
 since for each $\lambda \in \lau \Om {*-3} {\per p -\per e} {E_S/\sbat}^{\zdos} $ we have 
 $
  \zeta \wedge \lambda \in \lau \BOm{*-2} {\per{p}}{E_S}^{-\zdos}$ and $i_Z(\zeta \wedge \lambda ) = \lambda$.
 
 Consider now a good covering $\mathcal U$ of $S^{\sbat}$ and $\{f_U \tq U\in \mathcal U\}$ a subordinated partition of unity. The family $\{\tau_S^{-1}(U) \tq u \in \mathcal U\}$ is an open covering of $E_S$ and 
 $\{\tau_S \circ f_U \tq U\in \mathcal U\}$ is a subordinated partition of unity. These maps are $N$-invariant smooth maps constant on the fibers of $\tau_S$. This last property implies that $||f_U \circ \tau_S ||_S =||d(f_U \circ \tau_S )||_S=0$. So, the covering $\mathcal U$ possesses a subordinated partition of unity living in $\lau\BOm*{\per 0} {E_S}$.
 
 Using Bredon's trick \cite[p.~289]{MR1700700} one reduces the problem to the case $E_S = \R^{\dim S^{\sbat}} \times \R^{2b+2}$, where $\tau_S$ becomes the projection on the first factor. The action of the group $\sbat$ is trivial on the first factor. 
 
 Contracting this factor to a point, we reduce the problem to the case $E_S= \R^{2b+2} = \rc S^{2b+1} = (S^{2b+1}\times [0,\infty))/(S^{2b+1} \times \{0\})$ as filtered space. 
 Denoting by 
 $P$  the apex of the cone,
  we have $\per \chi (P) = 1, \per e(P)=2$ and  $\per p(P) =\per p( S) = p$.
 We need to prove that the inclusion
 $$
 I \colon \lau \BOm * {\per p - \per e} {\rc S^{2b+1}/\sbat}^{\zdos} \hookrightarrow
 \lau I * {\per p } {\rc S^{2b+1}} =
 \left\{i_Z \beta \tq   \beta \in \lau \BOm{*+1} {\per{p}}{\rc S^{2b+1}}^{-\zdos}\right\}
 $$
 is a quasi-isomorphism. This comes directly from \eqref{ayuda} with the equality $j_*Z =-Z$ (cf. \secref{twist} (b)).

\smallskip

 $\bullet$ {\em Step 3:  $\stra_3 =  \emptyset$}. Consider the invariant open covering $\mathcal V =\{T_S \tq S \in \stra_1\} \sqcup \{M\menos F_1\}$ of $M$. 
 We fix a smooth map $\lambda \colon [0,\infty) \to [0,1]$ verifying $\lambda=1$ on $[0,2]$ and $\lambda=0$ on $[3,\infty)$. 
 The map $f_S\colon T_S \to [0,\infty)$ is defined by $f_S(x) = \lambda(
 \rho_S(x))$ (see \subsecref{TMs}). It is an $\stres$-invariant smooth map, constant on the fibers of $\tau_S \colon D_S \to S$, which gives  $||f_S||_S = ||df_S||_S=0$. So, the family $\{ f_S \tq S \in \stra_1\} \sqcup \{ 1 - \sum f_S\}$ is a partition of unity, subordinated to $\mathcal V$, living in $\lau\BOm*{\per 0} M$.
  Now, it suffices to apply  Bredon's trick \cite[p.~289]{MR1700700} and the previous cases.

\smallskip

 $\bullet$ {\em Step 4:  $M = T_Q$ for some $Q \in \stra_3$}. Recall that $\tau_Q \colon T_Q \to Q$ is an $\stres$-invariant smooth bundle whose fiber is $\R^{f+1}$ for some $f \geq 3$. In fact,  the group $\stres$ acts trivially on $Q$ and orthogonally on the fiber $\R^{f+1}$ having the origin as the only fixed point. Notice that the action of $\stres$ on the sphere $S^f$ is a mobile action.

 Consider now a good covering $\mathcal U$ of $Q$ and $\{f_U \tq U\in \mathcal U\}$ a subordinated partition of unity. The family $\{\tau_Q^{-1}(U),\tq u \in \mathcal U\}$ is an open covering of $T_Q$ having 
 $\{f_U \circ \tau_Q\tq U\in \mathcal U\}$ a subordinated partition of unity. These maps are $\stres$-invariant smooth maps constant on the fibers of $\tau_Q$. This last property implies that $||f_U \circ \tau_Q||_Q =||d(f_U \circ \tau_Q)||_Q=0$ (cf. \eqref{TMe}).
 
 Using  Bredon's trick \cite[p.~289]{MR1700700} one reduces the problem to the case $M = \R^{\dim Q} \times \R^{f+1}$. The action of the group $\stres$ is trivial on the first factor. Contracting this factor to a point, we reduce the problem to the case $M= \R^{f+1} = \rc S^f = (S^f \times [0,\infty))/(S^f \times \{0\})$.
 Here, the stratum
 $Q $ corresponds to the apex $\tv$ of the cone. We have $\per \chii (\tv) = 3$ and $\per e(v)=4$. 
 The number $p \in \per \Z$ is defined by $\per p(Q) =p$. Take $\per{p}(\tv)=p$.
 We need to prove that the inclusion
 $$
 I \colon \lau \BOm * {\per p - \per e} {\rc S^f/\stres} \hookrightarrow 
  \lau I* {\per p } {\rc S^f}
 =
\left\{ i_{X_3} i_{X_2}i_{X_1} \omega  \tq \omega \in \lau \BOm {*+3} {\per p} {\rc S^f}\right\}
$$
is a quasi-isomorphism. 
 Notice that $S^3$ just acts on the sphere $S^f$. We have:
 $$
  \begin{array}{lcl}
 \lau \Om {*<p-4 } {\per p - \per e} {\rc S^{f}/\stres} &=& \lau \Om {*<p-4} {\per p - \per e}{ \left(S^{f} /\stres\right) \times (0,\infty)}, \\[,2cm]

 \lau \Om {p -4}  {\per p - \per e} {\rc S^{f}/\stres} &=& \left\{ \beta \in  \lau \Om {p-4} {\per p - \per e} { \left(S^{f} /\stres\right) \times (0,\infty)} \tq d\beta \equiv 0 \hbox{ on }  S^{f} \times (0,2)\right\}, \\[,2cm]
 
\lau \Om {*> p -4}  {\per p - \per e} {\rc S^{f}/\stres} &=&  \lau \Om {*>p-4} {\per p - \per e} { \left(S^{f} /\stres\right) \times (0,\infty), \left(S^{f} /\stres\right) \times (0,2)}, \\[,2cm]

 \lau I {* <p-3}{\per p}{\rc S^{f}}&=& \left\{   i_{X_3} i_{X_2}i_{X_1}\omega \tq \omega \in \lau  \BOm {* + 3 <p}{\per p } {S^{f}\times  (0,\infty)}\right\}, \\[,2cm]
 
 \lau I {p-3}{\per p }{\rc S^{f}}  &=& \left\{   i_{X_3} i_{X_2}i_{X_1} \omega \tq \omega \in \lau  \BOm {p}{\per p } {S^{f}\times  (0,\infty)}  \tq d\omega \equiv 0 \hbox{ on }  S^{f}  \times (0,2)\right\}, \hbox{ and} \\[,2cm]
 
 \lau I {*> p-3}{\per p }{\rc S^{f}}  &=& \left\{   i_{X_3} i_{X_2}i_{X_1}\omega \tq \omega \in \lau  \BOm {*+3>p}{\per p } {S^{f}\times  (0,\infty),S^{f}\times  (0,2)} \right\}.
 \end{array}
 $$
 Since  $\lau \Om {p -4}  {\per p - \per e} {\rc S^f/\stres} \cap d^{-1}(0) =  \lau \Om {p-4} {\per p - \per e} { \left(S^f /\stres\right) \times (0,\infty)} \cap d^{-1}(0)
 $ then  Step 3 gives that
 $
 I^* \colon \lau H {*} {\per p - \per e} {\rc S^f/\stres} \to \coho H {*}{ \lau I \cdot {\per p} {\rc S^f} } $ is an isomorphism for  \fbox{$* < p-3$}.
 
\begin{itemize}
\item[] \fbox{$*=p-3$} Since $\lau H {p-3} {\per p - \per e} {\rc S^f/\stres} =0$ (cf. \cite[Proposition~3.1.1]{MR2210257}) then we need to prove $\coho H {p-3}{ \lau I * {\per p} {\rc S^f} }=0$, that is:
$$
  \left\{
  \begin{array}{l}
  \omega \in \lau\BOm {p} {\per p}  {S^{f} \times (0,\infty)} \\[,2cm]
  \hbox{ with } d\omega \equiv 0 \ \text{on } S^{f} \times (0,2) \text{ and } d i_X \omega =0
  \end{array}
  \right.
  \Longrightarrow
    \left\{
  \begin{array}{l} \exists \eta \in \lau \BOm {p-1} {\per p} {S^{f} \times (0,\infty)} \\[,2cm]
   \hbox{with }  i_X  \omega = d i_X  \eta.
    \end{array}
  \right.
  $$
If $p=0$ then we can consider $\eta=0$ since $\omega $ is constant. Let us  suppose $p\geq 1$. Consider $\eta' = \int_1^- \omega$. It is an element of $\coho \BOm {p-1}  {S^f \times (0,\infty)} $ since the action of $\stres $ on the $(0,\infty)$-factor is trivial.
   A straightforward calculation gives 
   \begin{equation}\label{omega1}
     \omega = \pr^* \omega(1) + d \eta'  + \int_1^- d\omega. 
   \end{equation}
  Here, $\omega(1)$ is the restriction of $\omega$ to $S^f \times \{1 \}$ and $\pr \colon S^f \times (0,\infty) \to S^f \times \{ 1\}$ is the map defined by $\pr(x,t) =(x,1)$.
  Applying that both $\omega$ and $\eta'$ are $S^3$-invariant and $L_X=i_Xd + di_X$, we have 
  \begin{eqnarray*}
  i_{X_3} i_{X_2}i_{X_1} \omega& =& i_{X_3} i_{X_2}i_{X_1} \pr^* \omega(1) + i_{X_3} i_{X_2}i_{X_1} d \eta'+ i_{X_3} i_{X_2}i_{X_1} \int_1^-d\om \\ &= & i_{X_3} i_{X_2}i_{X_1} \pr^*\omega(1) -  d i_{X_3} i_{X_2}i_{X_1} \eta' - \int_1^- di_{X_3} i_{X_2}i_{X_1} \om =   i_{X_3} i_{X_2}i_{X_1} \pr^*\omega(1) -  d i_{X_3} i_{X_2}i_{X_1} \eta' .
  \end{eqnarray*}
  By hypothesis the differential form $\pr^*\omega(1)$ is a cycle of $\lau \BOm p {\per p} {S^f \times (0,\infty)}$.  Condition $  \per p \leq \per t$ implies $p  \leq \codim Q - 2 = f- 1$. This gives the existence of $\eta'' \in \coho \BOm {p-1} {S^f \times (0,\infty)}$ with $\pr^*\omega(1) = d \eta''$. We finish the proof taking $\eta = -\eta' - \eta''$.

\item[] \fbox{$*\geq p-2$}. Since $\lau H {*\geq p-2} {\per p - \per e} {\rc S^f/\stres} =0$ (cf. \cite{MR2210257})  then we need to prove $\coho H {*\geq p-2}{ \lau I* {\per p} {\rc S^f} }=0$, that is:

\smallskip

\scalebox{.85}{$
  \left\{
  \begin{array}{l}
  \omega \in \lau\BOm {*+3\geq p+1} {\per p}  {S^{f} \times (0,\infty)} \hbox{ with } d i_{X_3} i_{X_2}i_{X_1}  \omega =0
 \\[,2cm]
  \hbox{and }  \omega \equiv 0 \ \hbox{on } S^{f} \times (0,2)   \end{array}
  \right.
  \Longrightarrow
    \left\{
  \begin{array}{l} \exists \eta \in \lau \BOm {*+2 \geq p} {\per p} {S^{f} \times (0,\infty)}   \hbox{ with }  i_{X_3} i_{X_2}i_{X_1}   \omega = d i_{X_3} i_{X_2}i_{X_1}  \eta \\[,2cm]
  \hbox{ and }   d\eta \equiv 0 \hbox{ on } S^{f} \times (0,2) .    \end{array}
  \right.
  $}
  
  \smallskip
 \noindent Same proof as before with $\omega(1)=0$.
\end{itemize}

\smallskip

 $\bullet$ {\em Final Step}. Consider the invariant open covering $\mathcal V =\{T_Q \tq Q \in \stra_3\} \sqcup \{M\menos F_3\}$ of $M$. 
 We fix a smooth map $\lambda \colon [0,\infty) \to [0,1]$ verifying $\lambda=1$ on $[0,2]$ and $\lambda=0$ on $[3,\infty)$. 
 The map $f_Q\colon T_Q \to [0,\infty)$ is defined by $f_Q(x) = \lambda(
 \rho_Q(x))$ (see \subsecref{TMs}). It is an $\stres$-invariant smooth map, constant on the fibers of $\tau_Q \colon D_Q \to Q$, which gives  $||f_Q||_S = ||df_Q||_S=0$ for each singular stratum $S$ (cf. \eqref{TMe}). So, the family $\{ f_Q \tq Q \in \stra_3\} \sqcup \{ 1 - \sum f_Q\}$ is a partition of unity, subordinated to $\mathcal V$, living in $\lau\BOm*{\per 0} M$.
Now, it suffices to apply  Bredon's trick \cite[p.~289]{MR1700700} and the previous cases.
 \end{proof}

 Let us study the complex $ \lau K* {\per p} M$. Since the complex $\lau \Om *{\per p} {M/\stres}$ is included in $\ker \fint$ we have the short exact sequence
 $$
 0
 \to \lau \Om * {\per p } {M/\stres}  \hookrightarrow  \lau K* {\per p} M  \to \frac{ \lau K* {\per p} M } {\lau \Om * {\per p} {M/\stres}}  \to 0.
 $$
 This allows us to compute the cohomology of $ \lau K* {\per p} M$ in terms of the intersection cohomology $\lau H*{\per p} {M/\stres}$ and the cohomology of the complex $\lau C * {\per p} M =  \displaystyle \frac{ \lau K* {\per p} M } {\lau \Om * {\per p} {M/\stres}}$. 
 In fact, we are going to prove that this cohomology is residual relatively to the strata of $ \stra_1$. The calculation of 
 $\coho H * {\lau C \cdot{\per p} M}$ is carried out in several steps through the following restrictions
 \be\label{reduc}
 \xymatrix{
 M \ar@{~>}[r] & T_Q \sqcup T_S \ar@{~>}[r]  & Q \sqcup T_S \ar@{~>}[r]  & Q \sqcup S^{\sbat}=  \per {S}^{\sbat},
 }
 \ee
 where $S$ ranges over the strata of $\stra_1$ and $Q \in \stra_3$ with $\per Q \subset S$.
 We proceed in several steps.
  
   \begin{lemma}\label{lema2}
 Let $\Phi \colon \stres \times M \to M$ be a mobile action with $\stra_3 =\emptyset$. 
 For any perversity $ \per p $ on $M$ with $\per 0 \leq \per p \leq \per t$ we have
 \be\label{igual}
 \coho H* { \lau C \cdot {\per p} M } =
\bigoplus_{S \in \stra_1}\lau H {*-2p_S-2} {} { S^{\sbat}}^{-(-1)^{p_S }\zdos},
 \ee
 where
  the number $p_S$ denotes the integer part of $\per{p}(S)/2$.
  \end{lemma}
 \begin{proof}
 We proceed in several steps.
 
 \smallskip
 
 $\bullet$ {\em Step 1: $\stra_1=\emptyset$}. We do not have any singular stratum.  We need to prove that the LHS of the equality \eqref{igual} is 0.
 Let $\langle \omega \rangle$ be a cycle of $\lau C * {\per p} M$.
  The canonical decomposition of $\omega$ is 
 $$
\om =\omega_0 +  
   \chii_{1} \wedge  \omega_1 +    \chii_{2} \wedge  \omega_2 +    \chii_{3} \wedge  \omega_3 +
 \chii_{1} \wedge  \chii_2 \wedge  \omega_{12} +     \chii_{1} \wedge  \chii_3 \wedge \omega_{13}+ 
  \chii_{2} \wedge  \chii_3 \wedge \omega_{23}
$$
(cf. \eqref{hamalau}).
The differential form $\eta =-\chii_1\wedge \omega_{23} + \chii_2\wedge \omega_{13} - \chii_3 \wedge \omega_{12}$ belongs to $ \lau K* {\per p} M$ since $\eta \in \coho \BOm *{M\menos \Sigma}$ (cf. \eqref{lll} and \eqref{invll}), and
  $i_{X_3} i_{X_2}i_{X_1} \eta =0$.

The differential form $\omega' = \omega -d\eta$ verifies $\omega'_{12} =\omega'_{13} =\omega'_{23} =\omega'_{123} =0$. Since $d \omega' = d\omega  \in \lau \Om*{\per p} {M/\stres}$ then $\omega'_1=\omega'_2 =\omega'_3=0$ and therefore $\omega' =\omega'_0 \in \lau\Om *{\per p}{M/\stres}$. So, $ [\langle\omega\rangle] = [\langle\omega'\rangle ]= [0]=0$.

 \smallskip
 
 $\bullet$ {\em Step 2: $M=T_S$ for some $S\in \stra_1$}. In this case $\stra_1 = \{S\}$ and $S$ is a closed stratum. Define $E_S$ as in step 2 of the proof of Proposition \ref{31}.
 Using \eqref{intedef} and Corollary \ref{simpli} we get
 \begin{eqnarray*}
  \lau K* {\per p} M &=&
 \lau \BOm* {\per{p}}{E_S}^{\zdos}
 \oplus \left\{  \beta \in  \lau \BOm{*-2} {\per{p}}{E_S}^{-\zdos} 
  \tq i_Z\beta =0\right\} 
\oplus
 \left\{ \xi \in 
\lau \Om{*-1} {\per{p}}{E_S}^{-\zdos} \tq \ib L {Z} \ib L {Z} \xi=-\xi  
\right\}
\\
&\stackrel{\eqref{decompeseuno}} =&
 \lau \BOm* {\per{p}}{E_S}^{\zdos}
 \oplus  \lau \Om{*-2} {\per{p}}{E_S/\sbat}^{-\zdos} 
\oplus
 \left\{ \xi \in 
\lau \Om{*-1} {\per{p}}{E_S}^{-\zdos} \tq \ib L {Z} \ib L {Z} \xi=-\xi  
\right\},
\end{eqnarray*}
where the third complex is acyclic (see Corollary \ref{simpli}). So, using that  $ \lau \Om * {\per p} {M/\stres} \stackrel {\eqref{idde2}} =\lau \Om{*} {\per{p}}{E_S/\sbat}^{\zdos}$, we have that the cohomology of $\lau {C}*{\per p}{M}$ is that of the complex 
$$
\frac{\lau \BOm* {\per{p}}{E_S}^{\zdos}
 \oplus  \lau \BOm{*-2} {\per{p}}{E_S/\sbat}^{-\zdos} }{\lau \Om{*} {\per{p}}{E_S/\sbat}^{\zdos}},
 $$
\noindent relatively to the differential
 $
 D''\langle(\alpha,\beta)\rangle = \langle(d\alpha,d\beta - i_Z\alpha)\rangle . $
 Following \remref{rem}, this complex is quasi-isomorphic to 
 $$
 \frac{\lau \Om * {\per p} {E_S/\sbat}^{\zdos} \oplus \lau \Om {*-1} {\per p - \per e} {E_S/\sbat}^{-\zdos}
 \oplus  \lau \Om{*-2} {\per{p}}{E_S/\sbat}^{-\zdos} }{\lau \Om{*} {\per{p}}{E_S/\sbat}^{\zdos}}
 =
   \lau \Om {*-1} {\per p - \per e} {E_S/\sbat}^{-\zdos}
 \oplus  \lau \Om{*-2} {\per{p}}{E_S/\sbat}^{-\zdos} 
 $$
 endowed with the differential $D''(\lambda,\beta) = (-d\lambda, d\beta - \lambda)$.
It is immediately checked that the map $(\lambda,\beta) \mapsto  \langle \beta \rangle$ is a quasi-isomorphism between this complex and 
$
\lau \Om{*-2}{\per p/(\per p -\per e)}{E_S/\sbat}^{-\zdos}.
$
 The induced action $\Theta \colon N \times E_S \to E_S$ and the induced perversity $\per p$ on $E_S$  verify the conditions of \corref{lema}, with $\Sigma= S^{\sbat}$ and $q=p_S$. So, we get   
 $
 \lau H{*-2}{\per p/(\per p -\per e)}{E_S/\sbat}^{-\zdos} =\coho H {*-2p_S-2} {S^{\sbat}}^{-(-1)^{p_S }\zdos}.
 $

\smallskip

 $\bullet$ {\em Final step}. Just use Mayer-Vietoris as in  the Step 3 of the proof of the \propref{31}.
 \end{proof}

The first step of \eqref{reduc} comes from the following result concerning the Thom-Mather system (see \secref{TMs}):
 
  \begin{lemma}\label{lema22}
 Let $\Phi \colon \stres \times M \to M$ be a mobile action. 
 For any perversity $ \per p $ on $M$ with $\per 0 \leq \per p \leq \per t$ the restriction of forms 
 induces the quasi-isomorphism 
 $
 \lau C * {\per p} M \to \lau C * {\per p} {T_1 \cup T_3}.
 $ 
  \end{lemma}
 \begin{proof}
 Consider the invariant open covering $\{U = T_1 \cup T_3, V = M \menos \left( \rho_1^{-1}([0,3/4]) \cup \rho_3^{-1}([0,3/4])\right)\}$ of $M$. We have the Mayer-Vietoris sequence:
 $$
 \xymatrix{
 0 \ar[r] & \lau C * {\per p} {M} \ar[r] & 
 \lau C * {\per p} {U} \oplus \lau C * {\per p} {V}
 \ar[r] & 
 \lau C * {\per p} {U \cap V}
 \ar[r] & 0
 }.
 $$
 Following \lemref{lema2} we know that the complexes  $ \lau C * {\per p} {V}$ and $ \lau C * {\per p} {U \cap V}$ are acyclic. So, the restriction $ \lau C * {\per p} {M} \to 
 \lau C * {\per p} {U} $ is a quasi-isomorphism.
 \end{proof}

The term $Q \sqcup T_S$ appearing in \eqref{reduc} is not a manifold, but it is possible to define the complex $\lau C*{\per p}{-}$ on it using the following notion.

\begin{definition}\label{def35}
Let $\Phi \colon \stres \times M \to M$ be a mobile action. We consider a perversity $\per p$ on $M$. For each open subset $U \subset M$ we define
 $$
 \lau \Xi * {\per p} {U} =
 \{
 \omega \in \lau \BOm * {\per p} {U \menos \Sigma_3}  \tq \omega \ \hbox{ and } \ d \omega \hbox{ verify condition }\eqref{izarra}\},
 $$ 
 where this condition is
 \be\label{izarra}
  \om (v_0, \ldots , v_{\per p(Q)},-) = 0  
\hbox{ where $v_0, \ldots, v_{\per p(Q)}$ are  vectors tangent to the fibers of $ \tau_Q \colon( D_Q  \cap U )\to Q$},
 \ee
 for each $Q \in \stra_3$.
 We analogously define $ \lau \Xi * {\per p} {U/\stres} $ if $U \subset M$ is an $\stres$-invariant open subset.
 
 We define
 $$
 \lau {\widehat C}*{\per p} {T_S} =  \frac{ \lau {\widehat K}* {\per p} {T_S} = \lau K* {\per p} {T_S}\cap \lau \Xi *{\per p} {T_S}  }{\lau \Xi *{\per p} {T_S/\stres}}  .
$$

\end{definition}

We clearly have $\lau \Xi * {\per p} {U} = \lau \BOm * {\per p} {U}$ and 
$\lau {\widehat C} * {\per p} {U} = \lau C * {\per p} {U}$ if $U \supset T_3$ or if $\stra_3=\emptyset$.
 In particular, we have
 \be\label{KKgorro}
 \lau C  * {\per p} {T_1 \cup T_3}
 =
 \lau  {\widehat C} * {\per p} {T_1 \cup T_3}.
 \ee

  \begin{lemma}\label{lema3}
 Let $\Phi \colon \stres \times M \to M$ be a mobile action. 
 For any perversity $ \per p $ on $M$ with $\per 0 \leq \per p \leq \per t$ the restriction of forms
 induces the quasi-isomorphism 
 $
 \lau {\widehat C} * {\per p} {T_3 \cup T_1}\to 
 \lau {\widehat C} * {\per p} { T_1}. $ 
  \end{lemma}
 \begin{proof}
 Using Mayer-Vietoris, it suffices to prove that the restriction 
 $
 \lau {\widehat C} * {\per p} {T_3 }\to 
 \lau {\widehat C} * {\per p} { T_3 \cap T_1} $ 
 is a  quasi-isomorphism.
 Proceeding as in the Step 4 of the proof of \propref{31} we can suppose that  $T_3 
= \rc S^f$, where  $\stres$ acts orthogonally without fixed points on the sphere $S^{f}$, and trivially on the radius of the cone. The fixed points subset $F_3$ is, now, the apex $\tv$ of the cone. Moreover, in this new setting $T_3\cap T_1$ becomes $T_{\Sigma_1}\times(0,\infty)\subseteq\rc S^f$, being $\Sigma_1$ the union of non-mobile strata of the action in $S^f$ and $T_{\Sigma_1}$ a tubular neighbourhood of $\Sigma_1$ in $S^f$. 
We then have
$$
\lau {\widehat C} * {\per p} {T_3 } = \lau {\widehat C} * {\per p} {  S^f \times (0,\infty)  }\ \text{ and }\
\lau {\widehat C} * {\per p} {T_3 \cap T_1} = \lau {\widehat C} * {\per p} {  T_{\Sigma_1}\times (0,\infty)}.
$$
For these complexes, condition \eqref{izarra} becomes the following: 
$
  \om (v_0, \ldots , v_{p},-) = 0$ if $v_0, \ldots, v_{ p}$ are vectors  tangent to the first factor of $ (S^f \menos \Sigma_1) \times (0,2) $ 
or $ (T_{\Sigma_1}\menos \Sigma_1) \times (0,2)  $, respectively, being $p=\per p(\tv)$. 
In order to prove that the restriction $ \lau {\widehat C} * {\per p} {  S^f \times (0,\infty)  } \to  \lau {\widehat C} * {\per p} { T_{\Sigma_1} \times  (0,\infty) }$ is a quasi-isomorphism we notice that 
\begin{equation}\label{truncando}
	\begin{array}{lcl}
		\lau {\widehat C} {*< p} {\per p} {  S^f \times (0,\infty)  }  &=& \lau {C} {*< p} {\per p} {  S^f \times (0,\infty) },\\[,2cm]
		
		\lau {\widehat C} {p} {\per p} {  S^f \times (0,\infty)  }  &=& \left\{ \langle \beta \rangle \in  \lau {C} {*< p} {\per p} {  S^f \times (0,\infty)  } \tq d\beta \equiv 0 \hbox{ on }  S^{f}  \times (0,2)\right\},\text{ and} \\[,2cm]
		
		\lau {\widehat C} {* >p} {\per p} {  S^f \times (0,\infty)  } &=&  \lau {C} {* >p} {\per p} {  S^f \times (0,\infty) , S^f \times (0,2) }, \text{ whose cohomology is 0.}\\[,2cm]
	\end{array}
\end{equation}
Let us consider the operator $\pr \colon \lau {C} {*} {\per p} {  S^f }  \to  \lau {C} {*} {\per p} {  S^f \times (0,\infty)  } $ induced by the canonical projection. We define
the {\em truncated complex}
$
\sigma_p  \lau {C} {*} {\per p} {  S^f } = \lau {C} {<p} {\per p} {  S^f } \oplus ( \lau {C} {p} {\per p} {  S^f } \cap d^{-1}(0) ).
$
By \eqref{truncando}, we have that $\pr(\sigma_p  \lau {C} {*} {\per p} {  S^f })\subseteq  \lau {\widehat C} {*} {\per p} {  S^f \times (0,\infty)}$. Moreover, following \eqref{omega1} we get $\langle\omega \rangle =\pr\langle  \omega(1) \rangle + d \langle  \int_-^1 \omega \rangle  +  \langle  \int_-^1 d\omega \rangle$, and thus\newline $\pr\colon\sigma_p  \lau {C} {*} {\per p} {  S^f }) \to \lau {\widehat C} {*} {\per p} {  S^f \times (0,\infty)}$ is a quasi-isomorphism.
In a similar way, we prove  that the operator 
 \be\label{ayudabis}
  \pr \colon
 \sigma_p  \lau {C} {*} {\per p} {  T_{\Sigma_1} }  \to  \lau {\widehat C} {*} {\per p} {   T_{\Sigma_1} \times (0,\infty)  },
\ee
induced by the canonical projection, is a quasi-isomorphism.
So, the remaining task becomes to prove that the restriction
$\sigma_p \lau {C} * {\per p} {  S^f } \to  \sigma_p \lau {\widehat C} * {\per p} {T_{\Sigma_1}}$
induces a quasi-isomorphism, which is granted by \lemref{lema2}.  \end{proof}

By the definition of the complex $\lau {\widehat C}*{\per p} -$ we have the equality
$
\lau {\widehat C} *{\per p} {T_1 }= \oplus_{S \in \stra_1} \lau {\widehat C} *{\per p} {T_S}.
$

\begin{lemma}\label{355}
 Let $\Phi \colon \stres \times M \to M$ be a mobile action. For any perversity $\per 0 \leq \per p \leq \per t$ on $M$ and any stratum $S \in \stra_1$ we have
 \bee
 \coho H* { \lau  {\widehat C} \cdot {\per p} {T_S} } =
\lau H {*-2p_S-2} {\per {P_S}} {\per S^{\sbat}}^{-(-1)^{p_S }\zdos},
 \eee
  where
  the number $p_S$ denotes the integer part of $\per{p}(S)/2$.
 \end{lemma}
 \begin{proof} Consider the family $\mathcal Q = \{ Q \in \stra_3\tq Q \subset \per S\}$.
Recall that  $\per S^{\sbat}$ (cf. \secref{ss1}) is a filtered space. 
The regular part of $\per S^{\sbat}$ is $S^{\sbat}$ and its singular strata are the elements $Q \in \mathcal Q$. So,
 $$
 \lau \Om * {\per {P_S}} {\per S^{\sbat}} = 
 \Set{ \alpha \in \coho \Om * {S^{\sbat}} |
 	 \alpha \hbox{ and } d \alpha \hbox{ verify condition }\eqref{izarra}  \hbox{ for } \per P_S(Q) = \per p(Q) -  2p_S-2, \hbox{ where } Q \in \mathcal Q}.
 $$
Let us define the operator 
$$J_S \colon  \lau \Om {*-2p_S-2} {\per {P_S}} {\per S^{\sbat}}^{-(-1)^{p_S}\zdos } \to \lau {\widehat C} *{\per p} {T_S}
 $$
 by 
\be\label{JJS}
 J_S(\alpha) = \left<\gamma_1 \wedge \gamma_2 \wedge \tau_S^*\alpha \wedge  e_S^{p_S}\right>,
 \ee
 where  $e_S \in \lau \Om 2 {\per 2} {E_S/\sbat}^{-\zdos}$ is the Euler form of the $\sbat$-action on $E_S$ relatively to an $N$-invariant metric on $T_S$ (cf. \secref{twist}).
 
 \smallskip

 $\bullet$ {\em Step 1: The operator $J_S$ is well defined.}
 Since  $\tau_S \colon T_S \to S^{\sbat}$ is an $\stres$-equivariant map then we have $\tau_S^* \alpha  \in \lau \BOm {*-2p_S-2} {} {E_S}^{-(-1)^{p_S }\zdos } $ for each $\alpha \in \lau \Om {*-2p_S-2} 
 {\per {P_S}}{\per S^{\sbat}}^{-(-1)^{p_S}\zdos } $ and therefore 
  $\tau_S^*\alpha \wedge  e_S^{p_S} \in \coho \BOm {*-2} {E_S}^{-\zdos}$. For the perverse degree, we have
 $$
 ||\tau_S^*\alpha \wedge  e_S^{p_S} ||_S \leq  || e_S^{p_S} ||_S \leq  2 p_S  \leq\per p(S),
 $$
 and similarly for the differential $d(\tau_S^*\alpha \wedge  e_S^{p_S}  ) =
 \tau^*_S d\alpha\wedge 
  e_S^{p_S}  $. We conclude that $ \tau^*_S \alpha \wedge e_S^{p_S}  \in \lau \BOm * {\per p}{E_S}^{-\zdos}$.
 Applying \corref{forsimpli} we get that $\gamma_2 \wedge \gamma_3 \wedge \tau^*_S \alpha \wedge e_S^{p_S}  \in \lau \BOm * {\per p}{T_S = \stres \times_N E_S}$ with $\fint \gamma_2 \wedge \gamma_3 \wedge e_S^{p_S} \wedge \tau^*_S \alpha=0$. 
 
 It remains to prove that the differential form 
 $\eta =\gamma_2 \wedge \gamma_3 \wedge \tau^*_S \alpha \wedge e_S^{p_S} $ and its differential
 $d\eta =\gamma_2 \wedge \gamma_3 \wedge \tau^*_S d\alpha \wedge e_S^{p_S} $
  verify condition \eqref{izarra}  for  $ \per p(Q)$, where $Q \in \mathcal Q$. Consider a family 
  $v_0, \ldots, v_{\per p(Q)}$ of vectors tangent to the fibers of $ \tau_Q \colon( D_Q  \cap (T_S\menos S) )\to Q$. 
 Up to a reordering we get that $\eta(v_0, \ldots, v_{\per p(Q)})$ is a multiple of $\alpha (\tau_{S*}(v_0), \ldots, \tau_{S*}(v_{\per p(Q)-2p_S-2}))$. Conditions \eqref{TMe} imply that the vectors $\tau_{S*}(v_\bullet)$ are tangent to the fibers of $ \tau_Q \colon D_Q  \cap S \to Q$. Since $\alpha$ verifies condition \eqref{izarra}  for  $ \per {P_S}(Q)$, then we get 
 $\alpha (\tau_{S*}(v_0), \ldots, \tau_{S*}(v_{\per p(Q)-2p_S-2 = \per{P_S}(Q)}))=0$. So, $\eta$  verifies condition \eqref{izarra}  for  $ \per p(Q)$. The same argument applies to $d\eta$.
 
  We conclude that $\gamma_2 \wedge \gamma_3 \wedge\tau^*_S \alpha \wedge e_S^{p_S} \in \ker\fint$. The operator $J_S$ is well defined.

  \smallskip

 $\bullet$ {\em Step 2: The operator $J_S$ is a quasi-isomorphism when $M=T_Q$, for any $Q \in \stra_3$ satisfying $Q \subset \per S$.}
 Proceeding as in the Step 4 of the proof of \propref{31} we can suppose that $M =T_Q = \rc S^f$  where  $\stres$ acts orthogonally without fixed points on the sphere $S^{f}$ and trivially on the radius of the cone. 
 Notice that the action of  $\stres$ on the sphere $S^f$ is a mobile action. The stratum
 $Q$ is the apex $\tv$ of the cone. We  have $S=\Sigma_0\times(0,\infty)\subseteq\rc S^f$, being $\Sigma_0\subseteq S^f$ the union of some semi-mobile strata of the $S^3$-action in $S^f$.
 
 We also have $T_S=T_{\Sigma_0}\times (0,\infty)$. Since  $\per S = \rc (\Sigma_0)$ and $ \per S^{\sbat} = \rc ({\Sigma_0}^{\sbat} )$, 
 then $J_S$ becomes 
 \begin{equation}\label{JS36}
 J_S \colon 
\lau \Om {*-2p_S-2} {\per {P_S}} {\rc ({\Sigma_0}^{\sbat})}^{-(-1)^{p_S}\zdos } \to \lau {\widehat C} *{\per p} { T_{\Sigma_0}\times (0,\infty)}. 	
 \end{equation}
Let's prove that \eqref{JS36} is a quasi-isomorphism. 
On one hand, we know from \cite[Proposition~3.1.1]{MR2210257} that the operator $\pr \colon \sigma_{q}  \lau \Om {*} {\per {P_S}} {{\Sigma_0}^{S^1}} \to  \lau \Om {*} {\per {P_S}} {\rc ({\Sigma_0}^{S^1})}$, induced by the canonical projection, is a quasi-isomorphism. Here, $q = \per {P_S} (Q) = \per p (\tv) - 2p_S -2 = p -2p_S -2$.  On the other hand, using \eqref{ayudabis} we conclude that it suffices to prove that
$$
 J_S \colon 
  \sigma_q\lau \Om {*-2p_S-2} {\per {P_S}} {{\Sigma_0}^{S^1}}^{-(-1)^{p_S}\zdos } \to \sigma_p\lau {C} *{\per p} { T_{\Sigma_0} }
 $$
 is a quasi-isomorphism, which is granted by \lemref{lema2}.

\smallskip

 $\bullet$ {\em Final step.} Consider the invariant open covering $\mathcal V =\left\{  T_S\cap T_3,  T_S \menos  \rho_3^{-1}([0,2])\right\}$ of $T_S$. 
 We fix a smooth map $\lambda \colon [0,\infty) \to [0,1]$ verifying $\lambda=1$ on $[0,3]$ and $\lambda=0$ on $[4,\infty)$. 
 The map $f\colon T_S \to [0,\infty)$ is defined by $f(x) = \lambda(
 \rho_3(x))$. It is an $\stres$-invariant smooth map, constant on the fibers of $\tau_3 \colon D_3 \to \Sigma_3$, which gives  $||f||_S = ||df||_S= ||f||_Q = ||df||_Q=0$ for each $Q \in \stra_3$ with $Q \subset \per S$. So, the family $\{ f,1 - f\}$ is a partition of unity, subordinated to $\mathcal V$, living in $\lau\BOm*{\per 0} {T_S}$.
Now, it suffices to apply  Bredon's trick \cite[p.~289]{MR1700700} and the previous cases.
 \end{proof}

 \begin{proposition}\label{36}
 Let $\Phi \colon \stres \times M \to M$ be a mobile action. For any perversity $\per 0 \leq \per p \leq \per t$ 
 on $M$ we have
 \bee\label{igualbis}
 \coho H* { \lau C * {\per p} M } = 
\bigoplus_{S \in \stra_1}\lau H {*-2p_S-2} {\per {P_S}} {\per S^{\sbat}}^{-(-1)^{p_S }\zdos}.
 \eee
 \end{proposition}
 \begin{proof} It suffices to consider \lemref{lema22}, \eqref{KKgorro}, \lemref{lema3}  and \lemref{355}.
   \end{proof}

  \section{Gysin braid for a mobile action}\label{mobile}
  
  We construct two Gysin sequences associated with a mobile action $\Phi \colon \stres \times M \to M$. These sequences establish a relationship between the cohomology of the manifold $M$ and the intersection cohomology of the orbit space $M/\stres$. The existence of two distinct approaches to the cohomology of $M$ by the intersection cohomology of $M/\stres$ gives rise to two separate sequences: one from the left, and one from the right.

  The first approach, the  one from the left,  uses the short exact sequence
  \be\label{cortagys}
  0 \to \lau \Om * {\per p} {M/\stres} \to   \lau \BOm * {\per p} {M} \to \frac{ \lau \BOm * {\per p} {M}}{ \lau \Om * {\per p} {M/\stres}} \to 0,
  \ee
 where $\per p$ is a perversity on $M$.   The quotient $\displaystyle \lau {G} * {\per p} M = \frac{ \lau \BOm * {\per p} {M}}{ \lau \Om * {\per p} {M/\stres}}$ is the \emph{Gysin term}.
 
 \smallskip
 
  \begin{theorem}
 Let $\Phi \colon \stres \times M \to M$ be a mobile action. For any perversity $\per 0 \leq \per p \leq \per t$ on $M$ we have the long exact sequence, known as a \emph{Gysin sequence},
\be\label{gys1}
\xymatrix{
\cdots \ar[r] & \coho H {*-1} M \ar[r]  &
 \coho H {*-1} { \lau {G} \cdot {\per p } M }
   \ar[r] &  \lau H {*} {\per p} {M/S^3}  
\ar[r] & \coho H {*} M \ar[r]  & \cdots .
&
}
\ee
The cohomology of the Gysin term $ \lau {G} * {\per p } M$ is determined by the long exact sequence
\be\label{gys2}
\xymatrix@C=4mm{
\cdots \ar[r] & \coho H {*-1} { \lau {G} \cdot {\per p } M } \ar[r]  &
 \lau H {*-4} {\per p - \per e} {M/S^3} 
   \ar[r] &  
 \displaystyle \bigoplus_{S \in \stra_1}
\displaystyle  \lau H {*-2p_S-2} {\per{P_S}} {\per S^{\sbat}}^{-(-1)^{p_S}\zdos}
\ar[r] &  \coho H {*} { \lau {G} \cdot {\per p } M } \ar[r]  & \cdots .
&
}
\ee
 \end{theorem}
 \begin{proof}
 The exact sequence \eqref{gys1} comes from \eqref{cortagys} and from the fact that the complex $\lau \BOm*{\per p} M$ computes the cohomology of $M$ since  $\per 0 \leq \per p \leq \per t$ (cf. \secref{13} and  \propref{inv}).
We now consider  the short exact sequence
 \be\label{otra}
 \xymatrix{
 0 \ar[r] & \lau C * {\per p} M \ar[r]   &
 \lau {G} * {\per p } M  \ar[r]^{\per \fint}  &\lau I {*-3} {\per p} M \ar[r] & 0,
} 
\ee
 where $\per \fint \langle \omega\rangle =  \fint \omega$.
 The exact sequence \eqref{gys2} comes from  the fact that the cohomology of $ \lau I * {\per p} M$ is  $\lau H {*} {\per p - \per e} {M/S^3} $ (cf. 
\propref{31}) and  from \propref{36}.
 \end{proof}

  The second approach, the one from the right,  uses the short exact sequence
  \be\label{cortagysbis}
  0 \to  \lau K* {\per p} M \to   \lau \BOm * {\per p} {M} \to  \lau I {*-3} {\per p} M \to 0.
  \ee
   By symmetry, we say that the complex $ \lau K* {\per p} M$ is the \emph{co-Gysin term} of the action. 

\smallskip

  \begin{theorem}
 Let $\Phi \colon \stres \times M \to M$ be a mobile action. For any perversity $\per 0 \leq \per p \leq \per t$ on $M$ we have the long exact sequence, known as a \emph{Gysin sequence},
\be\label{gys3}
\xymatrix{
\cdots \ar[r] & \coho H {*-1} M \ar[r]  &
 \lau H {*-4} {\per p - \per e} {M/\stres} 
   \ar[r] &  \lau H {*} {} { \lau K\cdot {\per p} M } 
\ar[r] & \coho H {*} M \ar[r]  & \cdots 
}
\ee
The cohomology of the co-Gysin term $ \lau {K} \cdot {\per p } M$ is determined by the long exact sequence
\be\label{gys4}
\scalebox{.92}{
\xymatrix@C=6mm{
\cdots \ar[r] & \coho H {*-1} {  \lau K* {\per p} M } \ar[r]  &
 \bigoplus_{S \in \stra_1}
\displaystyle  \lau H {*-3-2p_S} {\per{P_S}} {\per S^{\sbat}}^{-(-1)^{p_S}\zdos}
\ar[r] &  
 \lau H {*} {\per p } {M/S^3} 
   \ar[r] & \coho H {*} {   \lau K* {\per p} M } \ar[r]  &
   \cdots .
&
}
}
\ee
 \end{theorem}
 \begin{proof}
 The sequence \eqref{gys3} is derived from \eqref{cortagysbis} and the following two facts:
\begin{itemize}
\item The complex $\lau \BOm*{\per p} M$ computes the cohomology of $M$, since $\per 0 \leq \per p \leq \per t$ (cf. \secref{13} and \propref{inv}).
\item The cohomology of the complex $\lau I {} {\per p} M$ is $\lau H {} {\per p - \per e} {M/\stres}$ (cf. \propref{31}).
\end{itemize}

The exact sequence \eqref{gys4} is derived from the short exact sequence
\begin{equation}
\label{OmKC}
0 \to \lau \Om*{\per p} {M/\stres} \to \lau K* {\per p} M \to \lau C* {\per p} M \to 0
\end{equation}
and \propref{36}.
 \end{proof}
 
 The previous long exact sequences \eqref{gys1}, \eqref{gys2}, \eqref{gys3} and \eqref{gys4} can be arranged in an exact braid diagram.
 
 \begin{definition}
 Let us consider six chain complexes $A^*,B^*,C^*,D^*,E^*$ and $F^*$.
 A {\em braid} is a diagram of chain maps of the form

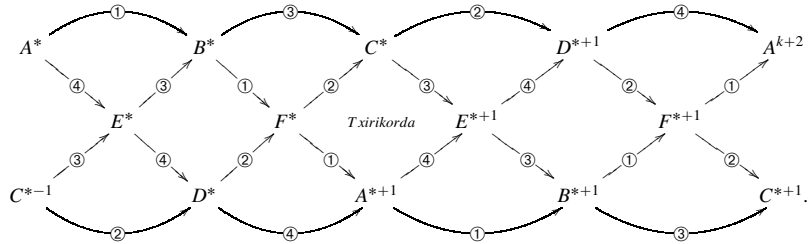
\bee
\xymatrix{
	A^*  \ar[rd] |\4\ar@/^{7mm}/[rr] |\1 &
	&
	B^*\ar[rd] |\1\ar@/^{7mm}/[rr] |\3&
	&
	C^* \ar[rd] |\3 \ar@/^{7mm}/[rr] |\2&
	&
	D^{*+1} \ar[rd] |\2 \ar@/^{7mm}/[rr] |\4&
	&
	A^{k+2} \\
	&
	E^* \ar[ur]  |\3\ar[dr] |\4 &
	&
	F^*\ar[ur]  |\2 \ar[dr] |\1  &
	&
	E^{*+1} \ar[ur] |\4\ar[dr] |\3 &
	&
	F^{*+1} \ar[ur] |\1\ar[dr] |\2&
	\\
	C^{*-1}   \ar[ur] |\3\ar@/_{7mm}/[rr] |\2&
	&
	D^* \ar[ur] |\2\ar@/_{7mm}/[rr]  |\4&
	&
	A^{*+1} \ar[ur] |\4\ar@/_{7mm}/[rr] |\1&
	&
	B^{*+1} \ar[ur] |\1\ar@/_{7mm}/[rr] |\3&
	&
	C^{*+1}.
}
\eee

 It is a {\em commutative braid} when all the triangles and diamonds are commutative. If the long sequences $\1$, $\2$, $\3$ and $\4$ are exact we say that the braid is an \em{exact braid}.
 \end{definition}
An exact and commutative braid possesses the two following properties.\label{propertiesB1B2}
\begin{enumerate}[B1-]
\item The following long sequence 
$$
\xymatrix@C=1cm{
\cdots \ar[r] & 
E^* \ar[r]^-{( \3 , \4)} &
 B^* \oplus D^*  \ar[r]^-{\1 - \2}&  F^*  \ar[r]^-{\3 \2} & E^{*+1}  \ar[r] & \cdots
}
$$
is exact (see for example \cite[pp.~39-41]{MR0178036}).
\item[]
\item The top and bottom sequences of the  braid are semi-exact sequences and both have the same exactness defaults: $\ker \3/\Ima\1 = \ker \4 / \Ima \2, \ldots$ (see for example \cite[p.~148]{MR755006}).
\end{enumerate} 

  \begin{theorem}\label{C}
 Let $\Phi \colon \stres \times M \to M$ be a mobile action. For any perversity $\per 0 \leq \per p \leq \per t$ on $M$ we have the following exact commutative braid, called the \emph{Gysin braid}:
 
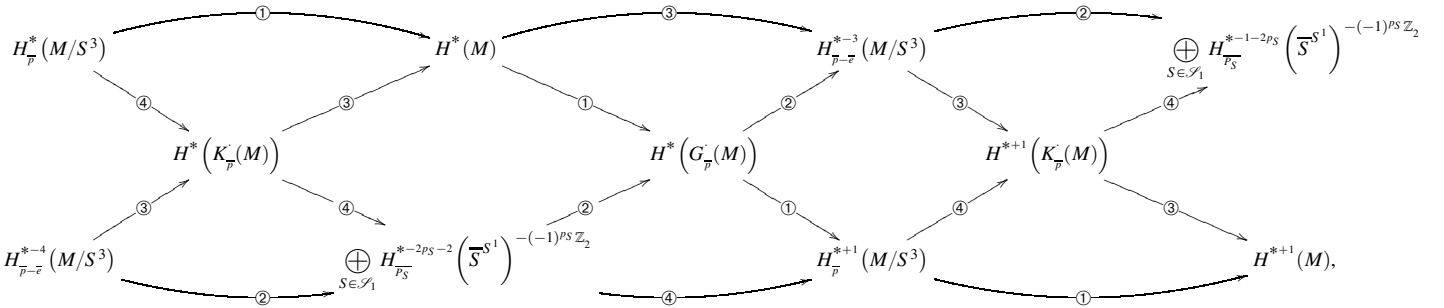
\bee
\scalebox{.7}{
 \xymatrix{
\lau H {*} {\per p} {M/\stres} \ar[rd] |\4\ar@/^{7mm}/[rr] |\1 &
&
\coho H{*}M \ar[rd] |\1\ar@/^{7mm}/[rr] |\3&
&
\lau H {*-3} {\per p - \per e} {M/\stres} \ar[rd] |\3 \ar@/^{7mm}/[rr] |\2&
&
\displaystyle \bigoplus_{S \in \stra_1} \lau H {*-1-2p_S} {\per{P_S}} {\per S^{\sbat}}^{-(-1)^{p_S}\zdos}  \\
&
\coho H * { \lau K \cdot {\per p} M}  \ar[ur]  |\3\ar[dr] |\4 &
&
\coho H {*} { \lau {G} \cdot {\per p } M } \ar[ur]  |\2 \ar[dr] |\1&
&
\coho H {*+1} { \lau K \cdot {\per p} M}   \ar[ur] |\4\ar[dr] |\3 &
\\
\lau H {*-4} {\per p - \per e} {M/\stres}  \ar[ur] |\3\ar@/_{7mm}/[rr] |\2&
&
\displaystyle \bigoplus_{S \in \stra_1} \lau H {*-2p_S-2} {\per{P_S}} {\per S^{\sbat}}^{-(-1)^{p_S}\zdos} \ar[ur] |\2\ar@/_{7mm}/[rr]  |\4&
&
\lau H {*+1}  {\per p} {M/\stres}\ar[ur] |\4\ar@/_{7mm}/[rr] |\1&
&
\coho H{*+1}M,
}}
\eee
\end{theorem}

\begin{proof}
In \cite{MR206943}, a braid is constructed from a triple. We follow this method to construct a braid associated with the following three complexes: $\lau \Om*{\per p} {M/\stres} \subset \lau K * {\per p} M \subset \lau \BOm*{\per p} {M}$. Recall that the cohomology of $\lau \BOm*{\per p} {M}$ is  $\lau H*{} {M}$ (cf. \secref{13}). To recognize the relative terms, we can do the following:
\begin{itemize}
\item The quotient $\frac{\lau \BOm*{\per p} {M}}{\lau \Om*{\per p} {M/\stres} }$ is the Gysin term $\lau {G} \cdot {\per p } M $.
   
\item    Using \eqref{cortagysbis} and \propref{31} we can determine the cohomology of $\frac{\lau \BOm*{\per p} {M}}{ \lau K * {\per p} M}$, which is given by $\lau H {*} {\per p - \per e} {M/\stres}$.

\item The cohomology of $\frac{ \lau K * {\per p} M}{\lau \Om*{\per p} {M/\stres} }$ is $\bigoplus_{S \in \stra_1} \lau H {*-2p_S-2} {\per{P_S}} {\per S^{\sbat}}^{-(-1)^{p_S}\zdos}$, as stated in \propref{36}.\qedhere
\end{itemize}
\end{proof}

\bigskip

 In some cases, the sequence \eqref{gys2}  splits  and the sequence \eqref{gys1} 
 is closer to the  classical  Gysin sequence. In particular, we  find the Gysin sequence \eqref{NuestraGysin} of  \cite{MR3119667}.

\begin{corollaire}\label{split}
Let $\Phi \colon \stres \times M \to M$ be a mobile action. When $\per p= \per 0$ the long exact sequence 
\eqref{gys2} splits   on the connecting homomorphism and we have the  
long exact sequence
\bee\label{bat}
\scalebox{.9}{
\xymatrix{
\cdots \ar[r] & \coho H {*-1} M \ar[r]  &
 \lau H {*-4} {} {M/S^3,\Sigma/\stres}  \oplus \displaystyle  \lau H {*-3} {} {M^{\sbat}}^{-\zdos}
   \ar[r] &  \lau H {*} {} {M/S^3}  
\ar[r] & \coho H {*} M \ar[r]  & \cdots .
&
}
}
\eee
\end{corollaire}

\begin{proof}
Since  $\per p = \per 0$ then $p_S =0$ and $\per{P_S} = \per{-2}$ for each $S \in \stra_1$.
 We have,
 \begin{eqnarray}\label{exoticoM}
  \bigoplus_{S \in \stra_1} \lau H {*-2p_S-3}  {\per {P_S}} {\per S^{\sbat}}^{-(-1)^{p_S} \zdos}
 &=&
 \bigoplus_{S \in \stra_1} \lau H {*-3} {\per {-2}} {\per S^{\sbat}}^{- \zdos}
 \stackrel{\hbox{ \tiny\cite[Prop. 13.5 ]{CST5}}}{=}
  \bigoplus_{S \in \stra_1} \lau H {*-3} {} {\per S^{\sbat} , \per S^{\stres}}^{- \zdos}
 \\ \nonumber
 & = &
 \lau H {*-3} {} {M^{\sbat},M^{\stres}}^{-\zdos}
=
 \lau H {*-3} {} {M^{\sbat}}^{- \zdos}
 \end{eqnarray}
since $\coho H {*-3}{M^{\stres}}^{-\zdos}=0$. 

Consider a cycle $\omega \in \lau I * {\per 0} M$. Let's compute the connecting homomorphism $\delta[\omega] $ of the sequence \eqref{gys2}.
We have seen in the proof of \propref{31} that $\chi_3 \wedge \chi_2 \wedge \chi_1 \wedge \omega \in 
\lau \BOm {*+3} {\per p} M$ with $\fint (\chi_3 \wedge \chi_2 \wedge \chi_1 \wedge \omega) = \omega$.
Using \eqref{ddd} we get
$$
d(\chi_3 \wedge \chi_2 \wedge \chi_1) \wedge \omega = (e_1^2+e_2^2+e_3^2) \wedge \omega -
d((e_3 \wedge \chi_3 +e_2 \wedge \chi_2 +e_1 \wedge \chi_1)\wedge \omega ).
$$
Since $\omega \in \lau \Om*{\per 0 -\per \chi} {M/\stres}$ then $\omega$ vanishes on $D_1 \sqcup D_3$ and therefore, so does  any multiple of $\omega$. So, $(e_1^2+e_2^2+e_3^2) \wedge \omega \in \lau \Om*{\per 0 -\per \chi} {M/\stres}$ and
$ (e_3 \wedge \chi_3 +e_2 \wedge \chi_2 +e_1 \wedge \chi_1)\wedge \omega  \in \lau \BOm * {\per 0} M$. This gives 
\begin{eqnarray*}
d \langle \chi_3 \wedge \chi_2 \wedge \chi_1 \wedge \omega \rangle 
&=&
 \langle d(\chi_3 \wedge \chi_2 \wedge \chi_1) \wedge \omega \rangle 
 =
 \langle 
 (e_1^2+e_2^2+e_3^2) \wedge \omega -
d((e_3 \wedge \chi_3 +e_2 \wedge \chi_2 +e_1 \wedge \chi_1)\wedge \omega \rangle\\
&=&
- \langle 
d(e_3 \wedge \chi_3 +e_2 \wedge \chi_2 +e_1 \wedge \chi_1\wedge \omega )
\rangle,
\end{eqnarray*}
yielding $\delta[\omega] =0$.
Now, the sequence \eqref{bat} comes from \eqref{gys1} (cf. \cite[Proposition 13.4)]{CST5}), \eqref{otra} and \eqref{exoticoM}.
\end{proof}

  \begin{remarques}\label{Rem}
 
 \color{white}.\negro
  
  $(a)$ \emph{The exotic term 
 $\displaystyle \bigoplus_{S \in \stra_1} \lau H * {\per {P_S}} {\per S^{\sbat}}^{-(-1)^{p_S} \zdos}$ }. When the exotic terms vanishes the Gysin Braid becomes the long exact sequence \eqref{gysfinal}. This happens when $\stra_1 = \emptyset$. In particular,  when the action $\Phi$ is  almost-free (i.e., $\stra_1=\stra_3=\emptyset$) or semi-free.
 
We have seen in the proof of Corollary \ref{split} that the   exotic term can be simplified when $\per p = \per 0$.
 In this case we have  $
  \bigoplus_{S \in \stra_1} \lau H * {\per {P_S}} {\per S^{\sbat}}^{-(-1)^{p_S} \zdos} =
 \lau H * {} {M^{\sbat}}^{- \zdos}.
 $

  The isotropy subgroup of a point of $M^{\sbat}\menos M^{\stres}$ is conjugated to $\sbat$ or $N$. Let us suppose that the first situation does not appear, that is, the group $\zdos$ acts trivially on $M^{\sbat}$. 
  This implies that $\lau H * {} {M^{\sbat}}^{- \zdos}=0$.
  
  Equality \eqref{exoticoM} is  still true when both $\per p(S)\equiv 0,1\pmod{4}$ for all $S\in\stra_1$ and $\per p(S)\le 4p+1$ for all $S\in\stra_3$.

We have  
$$
\lau H * {} {M^{\sbat}}^{- \zdos} =
\lau H * {} {M^{\sbat}}/  \lau H * {} {M/\stres},
$$
since
$ \lau H * {} {M^{\sbat}}=  \lau H * {} {M^{\sbat}}^{ \zdos} \oplus  \lau H * {} {M^{\sbat}}^{- \zdos} =
 \lau H * {} {M^{\sbat}/\zdos} \oplus  \lau H * {} {M^{\sbat}}^{- \zdos} 
 =
  \lau H * {} {M/\stres} \oplus  \lau H * {} {M^{\sbat}}^{- \zdos} .
 $

Notice that $\displaystyle \bigoplus_{S \in \stra_1} \lau H * {\per {P_S}} {\per S^{\sbat}}^{-(-1)^{p_S} \zdos} =\bigoplus_{S \in \stra_1} \lau H * {} {S^{\sbat}}^{-(-1)^{p_S} \zdos}$ when $\stra_3=\emptyset$. 

 \smallskip

 $(b)$ Property B1 (cf. p.~\pageref{propertiesB1B2}) yields the following two long exact sequences:
 \bee
 \xymatrix@=.6cm{
 \cdots \ar[r] & \coho H {*} { \lau K \cdot {\per p} M}   \ar[r] & \coho H {*} { \lau G \cdot {\per p} M}   \ar[r] & \lau H {*-3} {\per p - \per e} {M/\stres}\oplus \lau H {*+1} {\per p} { M/\stres}   \ar[r] 
 &
 \coho H {*+1} { \lau K \cdot {\per p} M}   \ar[r] &
 \cdots \\
  \cdots \ar[r] & \coho H {*-1} { \lau G \cdot {\per p} M}   \ar[r] & \coho H {*} { \lau K \cdot {\per p} M}   \ar[r] & \lau H {*} {} {M}\oplus \displaystyle \bigoplus_{S \in \stra_1} \lau H {*-2p_S-2} {\per{P_S}} {\per S^{\sbat}}^{-(-1)^{p_S}\zdos}  \ar[r] 
 &
 \coho H {*} { \lau G \cdot {\per p} M}   \ar[r] &
 \cdots 
 }
 \eee
 relating the cohomologies of the Gysin and co-Gysin terms.
 
 \smallskip

 $(c)$ Consider the case $\per p=\per 0$. Since the bottom map \2  of the Gysin braid vanishes (cf. \corref{split}) a diagram chasing gives the short exact sequence
   $$
  \xymatrix@C=1.5cm{
 0 \ar[r] & \Ima \1_{top}  \ar@{^(->}[r] &  \ker \3_{top} 
    \ar[r]^-{\2^{-1} \circ \1} & \ker \4_{bottom} \ar[r] & 0,
   }
   $$
 (cf. B2). 

\smallskip

 $(d)$ The splitting property of Corollary \ref{simpli} may not hold if $\per p\ne\per 0$, as the following example shows.

 Let us consider the join  $M = S^a  * S^2 * S^3 = S^{a+7} $ with $a\geq 1$ (see \cite[p.~468]{MR1700700}). 
  The action $\Phi \colon \stres \times M \to M$ is defined in each factor as follows:
  \begin{itemize}
   \item $S^3$ acts trivially on $S^a$.
   \item  $\stres$ acts on the left of the homogeneous space $S^2 = \stres/\sbat$.
   \item $\stres$  acts on the left of $S^3$ by multiplication on $S^3$. 
   \end{itemize}
   We put $\{b_1,b_2\}$ the two points of $S^2$ whose isotropy subgroup is $\sbat$.
  
\smallskip

We have $\stra_3 = \{Q = S^a\}$ and $\stra_1 = \{ S=(S^a * S^2) \menos S^a\}$.
Notice that 
  $$
  M^{\sbat} ={ \per S}^{\sbat}  = \left(S^a * S^2 \right)^{\sbat} =  S^a * \left(S^2\right)^{\sbat}=  S^a * \{b_1,b_2\},
  $$
the action of  $g \in \zdos$ interchanges both points. 
  
  The orbit space $M/\stres$ is the filtered space $S^a * \left( S^2 * S^3\right)/\stres   =
   S^a *  \Sigma S^2 = S^{a+4}$  endowed with the filtration
     $S ^a \subset S^a * \{P\} \subset S^a * \Sigma S^2 = S^{a+4}$ where $P$ is one of the two apices of the suspension 
     $\Sigma S^2$.

The perversity $\per p$ is given by two numbers $(p_1,p_3) = (\per p(S),\per p (Q))$.
Since  $\dim M = a +7$, $\dim Q=a$ and $\dim S=a +3$ then  the condition $  \per 0 \leq \per p \leq \per t$ becomes $(0,0 )\leq (p_1,p_3) \leq (2,5)$. In particular, the perversities $\per 0 =(0,0)$ and $\per e = (2,4)$ satisfy this condition. If $\per p =\per 0$ (resp. $ \per e$ ) then $\per {P_S}= \per{-2}$ (resp. $\per{0}$).

 A straightforward calculation using  \cite[Proposition 13.5]{CST5} gives
 
 \smallskip
 
  \begin{tabular}{llll}
  $\bullet$  $H^i(M)$ & =& $\R $ & if $ i= 0,a+7 $,\\[,21cm]
    $\bullet$ 
  $\lau H i {-\per e} {M/\stres}  = \lau H i {} {S^{a+4} ,S^a * \{P\}}$  &=& $\R$ & if $i=a+4$,
  \\[,21cm]
    $\bullet$  $\lau H i {\per 0} {M/\stres}  = \coho H i{S^{a+4}}$  &=& $\R$ & if $i=0,a+4$ ,\\[,21cm]
  $\bullet$   $\lau H i { \per e} {M/\stres} = \lau H i {} {S^{a+4} \menos (S^a * \{P\})}$ &=& $\R$ &if $i=0$,
\\[,21cm]
     $\bullet$  $\lau H i {\per {0}}  {{\per S}^{\sbat}}^{-\zdos} 
  =\coho H i {S^{a+1} \menos S^a}^{-\zdos} $&=&$\R$ &
  if $i=0$, \text{and} \\[,21cm]
     $\bullet$  $\lau H i {\per {-2}}  {{\per S}^{\sbat}}^{\zdos} 
  =\coho H i {S^{a+1},S^a}^{\zdos} $&=&$\R$ & 
  if $i=a +1$,
  \end{tabular}
  \smallskip
  
\noindent 
and 0 for the other values of $i$.

One easily checks that the long exact sequence 
\eqref{gys2} does not split   on the connecting homomorphism when $\per p = \per e$. Also, this example shows that the long exact sequence  \eqref{gys4}  does not split   on the connecting homomorphism, even in the case  $\per p = \per 0$.

\smallskip

$(e)$ Let's recall other exact braids in the literature. For the classical case of a semifree action of $\stres$ both the Smith-Gysin and the Gysin sequence are part of an exact braid, along with the long exact sequences of the pairs $(M,M^{\stres})$ and $(M/\stres,M^{\stres}/\stres)$, for singular cohomology (see \cite[p.~161]{MR0413144}). This braid has recently been generalized in \cite{MR4731363} to the case of any $\stres$-action by using Verona forms (as in \cite{MR3119667}), which corresponds to the case ${\per p}=\per{0}$ of the present paper. When intersection cohomology is considered, we have seen that the Gysin sequence itself leads to an exact braid. It is reasonable to suspect that a similar phenomenon will happen with the Smith-Gysin sequence. The interplay between the braids is still to be explored.

\end{remarques}

 
 We end the Section by studying the behaviour of the Gysin and co-Gysin terms when the perversity changes.

  \begin{proposition}\label{D}
 Let $\Phi \colon \stres \times M \to M$ be a mobile action. For any perversities $\per 0 \leq \per p \leq \per q \leq  \per t$ on $M$ we have the exact commutative braids,
\bee
\scalebox{.8}{
 \xymatrix{
\lau H {*} {\per p} {M/\stres} \ar[rd] |\4\ar@/^{7mm}/[rr] |\1 &
&
\coho H*M \ar[rd] |\1\ar@/^{7mm}/[rr] |\3&
&
\lau H {*} { } { \lau G\cdot {\per q} M } \ar[rd] |\3 \ar@/^{7mm}/[rr] |\2&
&
\lau H {*+1} {\per q/ \per p} {M/\stres} 
 \\
&
  \lau H {*} {\per q} {M/\stres} \ar[ur]  |\3\ar[dr] |\4 &
&
\lau H {*} {} { \lau G\cdot {\per p} M }  \ar[ur]  |\2 \ar[dr] |\1&
&
\lau H {*+1} {\per q} {M/\stres}   \ar[ur] |\4\ar[dr] |\3 &
\\
\lau H {*-1} { } { \lau G\cdot {\per q} M }  \ar[ur] |\3\ar@/_{7mm}/[rr] |\2&
&
\lau H {*} {\per q/ \per p} {M/\stres}  \ar[ur] |\2\ar@/_{7mm}/[rr]  |\4&
&
\lau H {*+1}  {\per p} {M/\stres}\ar[ur] |\4\ar@/_{7mm}/[rr] |\1&
&
\coho H{*+1}M,
}}
\eee
and
\bee
\scalebox{.8}{
 \xymatrix{
\lau H {*} { } { \lau {K}\cdot {\per p} M }  \ar[rd] |\4\ar@/^{7mm}/[rr] |\1 &
&
\coho H*M \ar[rd] |\1\ar@/^{7mm}/[rr] |\3&
&
\lau H {*-3} {\per q -\per e } {M/\stres}  \ar[rd] |\3 \ar@/^{7mm}/[rr] |\2&
&
\lau H {*-3} {(\per q - \per e)/ (\per p - \per e)} {M/\stres} 
 \\
&
 \lau H {*} { } { \lau {K}\cdot {\per q} M }
  \ar[ur]  |\3\ar[dr] |\4 &
&
\lau H {*-3} {\per p -\per e } {M/\stres}    \ar[ur]  |\2 \ar[dr] |\1&
&
\lau H {*+1} { } { \lau {K}\cdot {\per q} M }  \ar[ur] |\4\ar[dr] |\3 &
\\
\lau H {*-4} {\per q - \per e} {M/\stres}  \ar[ur] |\3\ar@/_{7mm}/[rr] |\2&
&
\lau H {*-4} {(\per q - \per e)/ (\per p - \per e)} {M/\stres}  \ar[ur] |\2\ar@/_{7mm}/[rr]  |\4&
&
\lau H {*+1} { } { \lau {K}\cdot {\per p} M } \ar[ur] |\4\ar@/_{7mm}/[rr] |\1&
&
\coho H{*+1}M.
}}
\eee

\end{proposition}

\begin{proof}
	As in the proof of Theorem \ref{C}, we shall apply the technique of  \cite{MR206943} to construct two braids associated with the following pairs of triples: 
 $ \lau \Om*{\per p} {M/\stres} \subset  \lau \Om*{\per q} {M/\stres} \subset \lau \BOm * {\per q} M$ and
$\lau K* {\per p} M \subset  \lau K* {\per q} M\subset \lau \BOm * {\per q} M$. 
It is important to recall that the cohomology of $\lau \BOm*{\per p} {M}$ or $\lau \BOm*{\per q} {M}$ is $\coho H*M$ (cf. Section \ref{13}). 
To identify the relative terms, we consider the following schema
$$
\renewcommand{\arraystretch}{3.1}
\begin{array}{|c|c|c|c|c|c|} \hline
\frac{\lau \Om*{\per q} {M/\stres} }{\lau \Om*{\per p} {M/\stres}}  &\frac{\lau \Om*{\per q} {M} }{\lau \Om*{\per q} {M/\stres}}&
\frac{\lau \Om*{\per q} {M} }{\lau \Om*{\per p} {M/\stres}}&\frac{ \lau K* {\per q} M}{\lau K* {\per p} M  }&
\frac{\lau \BOm*{\per q} {M}}{ \lau K * {\per q} M}&
\frac{\lau \BOm*{\per q} {M} }{\lau K* {\per p} M }\\ \hline
 \lau \Om*{\per q/\per p} {M/\stres} &\lau G*{\per q} {M}&\lau G*{\per p} {M}&\lau \Om {*-4} {(\per q - \per e)/ (\per p - \per e)} {M/\stres} & \lau \Om {*} {\per q - \per e} {M/\stres}&
 \lau \Om {*} {\per p - \per e} {M/\stres} \\ \hline
 \end{array}
 $$
 The first two columns are actually equalities. In the other columns, we have complexes with the same cohomology.
To prove this fact for the third column, we use the equality $\lau \Om*{\per q}{M} \cap \lau \Om*{\per p}{M/\stres} = \lau \Om*{\per q}{M/\stres}$ and the following fact, proven in \cite{MR1143404,MR2210257}, which we shall refer to as the $(\per p,\per q)$ property: the inclusion $\lau \Om*{\per p}{M} \hookrightarrow \lau \Om*{\per q}{M}$ induces an isomorphism in cohomology.

For column 4, we use \eqref{cortagysbis} and \propref{31}. Column 5 comes from the equality $\lau K*{\per q}{M} \cap \lau \BOm*{\per p}{M} = \lau K*{\per p}{M}$ and the $(\per p,\per q)$ property. The last column comes from the short exact sequence
 $
 0\to \lau K {*} {\per p} {M}  \to \lau \Om* {\per p} M \to \lau I {*-3} {\per p}{M} \to 0,
 $
the $(\per p,\per q)$ property  and, finally, \propref{31}.
\end{proof}

\section{Gysin sequence for a non-mobile action}\label{non-mobile}
In this section, we consider a non-mobile non-trivial action $\Phi \colon \stres \times M \to M$. The family of regular strata 
(resp. singular strata) is $\stra_1\ne \emptyset$ (resp. $\stra_3$). A perversity is a map $\per p \colon \stra_3 \to \per \Z$, that is, a family of numbers $\{ \per p(Q)\tq Q \in \stra_3\} \subset \per \Z$.  We consider on $M/\stres$ and $M^{\sbat}$ the induced filtered space structure. Notice that the family of singular strata is still $\stra_3$.

  \begin{theorem}\label{ThD}
 Let $\Phi \colon \stres \times M \to M$ be a non-mobile non-trivial action. For any perversity $\per 0 \leq \per p \leq \per t$ on $M$ we have 
 \be\label{final}
\lau H * {} M = \lau H * {\per p}  {M/\stres}
\oplus
 \lau H {*-2} {\per p - \per 2}  {M^{\sbat} }^{-\zdos}. 
\ee
\end{theorem}
\begin{proof}
Since the action of $\stres$ on $M\menos \Sigma$ has no fixed points, the assignment $(g,x) \mapsto g \cdot x$, 
establishes an $\stres$-equivariant diffeomorphism between the twisted product $\stres \times_N \left( M^{\sbat}\menos \Sigma\right)$ and $M\menos \Sigma$. Here, $\Sigma = M^{\stres}$ is the union of singular strata.
Following \corref{forsimpli} we have 
$$
\lau \BOm * {} {M \menos \Sigma} = 
\lau \Om *{}{ M^{\sbat}\menos \Sigma}^{\zdos} \oplus \lau \Om {*-2}{}{ M^{\sbat}\menos \Sigma} ^{-\zdos},
$$
since $Z=0$. A differential form $\omega \in  \lau \BOm * {} {M \menos \Sigma} $ is $\alpha + \gamma_2 \wedge \gamma_3 \wedge \beta$ with $(\alpha,\beta) \in \lau \Om *{}{ M^{\sbat}\menos \Sigma}^{\zdos} \oplus \lau \Om {*-2}{}{ M^{\sbat}\menos \Sigma} ^{-\zdos}$. It remains to compute the perverse degree $||\omega||_Q$.

We consider an $\stres$-invariant Thom-Mather system $\mathfrak T_M = \{T_Q\tq Q \in \stra_3\}$. It induces the $N$-invariant Thom-Mather system 
$\mathfrak T_{M^{\sbat}} = \{T_Q \cap M^{\sbat} \tq Q \in \stra_3\}$ on $M^{\sbat}$ (cf. \secref{circuloN}).
Notice that $T_Q \menos Q = S^3 \times_N \left((T_Q \cap M^{\sbat})\menos Q\right)$.
The map $\tau_Q  \colon T_Q \menos Q \to Q$ becomes $\langle g,x \rangle \mapsto \tau_Q(x)$.
So, the fiber of $\tau_Q$ over a point $y \in Q$ is $S^3 \times_N \left( (\tau^{-1}_Q (y) \cap M^{\sbat}) \menos Q\right)$.
This gives
$
||\omega||_Q = \max \{ ||\alpha||_Q, 2 + ||\beta||_Q\}
$
and therefore 
$$  
\coho H * M = \coho H * {\lau \BOm \cdot {\per p} M}= \lau H * {\per p}  {M^{\sbat} }^{\zdos}
\oplus
 \lau H {*-2} {\per p - \per 2}  {M^{\sbat}  }^{-\zdos},
$$
since the complex $\lau \BOm*{\per p} M$ computes the cohomology of $M$ for the perversity  $\per 0 \leq \per p \leq \per t$ (cf. \secref{13} and  \propref{inv}).
Finally, we get  \eqref{final} from $(M^{\sbat}\menos \Sigma)/{\zdos} = (M\menos \Sigma) /\stres$.
\end{proof}
 
 \begin{remarque}
 Considering the perversity $\per p = \per 0$ we get
$$
\lau H * {} M = \lau H * {}  {M/\stres}
\oplus
 \lau H {*-2} {}  {M^{\sbat} }^{-\zdos}
$$
(cf. \eqref{exoticoM} and \cite[Proposition 13.4)]{CST5}).
\end{remarque}
 \bibliographystyle{amsplain} 

\bibliography{BaseS3}

\end{document}